\documentclass[11pt]{amsart}

\usepackage{amsmath,amssymb,amsthm}

\usepackage{amsmath, amssymb} 
\usepackage{amsthm, amsfonts, mathrsfs}
\usepackage{amsfonts,graphicx}

\theoremstyle{plain}
\newtheorem{Theo}{Theorem}[section]
\newtheorem{lem}[Theo]{Lemma}
\newtheorem{cor}[Theo]{Corollary}
\newtheorem{prop}[Theo]{Proposition}

\theoremstyle{plain}
\theoremstyle{definition}
\newtheorem{defi}[Theo]{Definition}

\theoremstyle{remark}
\newtheorem{Rema}[Theo]{Remark}
\newtheorem*{rema*}{Remark}
\newtheorem{remas}[Theo]{Remarks}

\newcommand{\ZZ}{\mathbb{Z}}  
\newcommand{\CC}{\mathbb{C}}  
\newcommand{\NN}{\mathbb{N}}

\newcommand{\RR}{\mathbb{R}}

\newcommand{\DD}{|\textnormal{D}|}

\author[T.   Hmidi]{Taoufik Hmidi}
\address{IRMAR, Universit\'e de Rennes 1\\ Campus de
Beaulieu\\ 35~042 Rennes cedex\\ France}
\email{thmidi@univ-rennes1.  fr}

\date{}

\begin{document}

\title[On the critical Boussinesq system]
{\small{On a maximum principle and its application to logarithmically critical Boussinesq system}}
\maketitle

\begin{abstract}
In  this paper we study a  transport-diffusion model with some  logarithmic dissipations. We look for two kinds of estimates. The first one is a maximum principle whose proof is based on Askey theorem concerning characteristic functions and some tools from the theory of $C_0$-semigroups. The second one is a smoothing effect based on  some results from harmonic analysis and sub-Markovian operators. As an application we prove the global well-posedness for the two-dimensional Euler-Boussinesq system where the dissipation occurs only on the temperature equation and has the form $\frac{\DD}{\log^\alpha(e^4+\DD)}$, with $\alpha\in[0,\frac12]$. This result improves the critical dissipation $(\alpha=0)$ needed for global well-posedness which was  discussed in  \cite{HKR}.  \end{abstract}
{\footnotesize{\tableofcontents}}
\maketitle

\section{Introduction} 
The first goal of this paper is to study some mathematical problems related to  the following transport-diffusion model with logarithmic dissipations
\begin{equation}
\label{trdif}
\left\{ 
\begin{array}{ll} 
\partial_{t}\theta+v\cdot\nabla\theta+ \kappa \frac{\DD^\beta}{\log^\alpha(\lambda+\DD)} \theta=0, (t,x)\in\RR_+\times\RR^d\\
\textnormal{div}\,v=0\\
 \theta_{| t=0}=\theta_{0}. 
\end{array} \right.
\end{equation}
Here, the unknown is the scalar function $\theta$,  the  velocity $v$ is a time-dependent vector field with zero divergence and $\theta_0$  is the initial datum. The parameter $\kappa\geq0$, $\lambda>1$ and $\alpha,\beta\in\RR.$ The operator $\frac{\DD^\beta}{\log^\alpha(\lambda+\DD)}$ is defined through its Fourier transform 
\centerline{
$
\mathcal{F}({\frac{\DD^\beta}{\log^\alpha(\lambda+\DD)}f})(\xi)=\frac{|\xi|^\beta}{\log^\alpha(\lambda+|\xi|)}(\mathcal{F}f)(\xi)
.$ }

 We will discuss along this paper some quantitative properties for this model, especially two kinds of information will be established: maximum principle  and  some  smoothing effects. We notice that the special case of the  equation \eqref{trdif} corresponding to  $\alpha=0$ and $\beta\in[0,2]$ appears naturally in some fluid  models like  quasi-geostrophic equations or  Boussinesq systems. In this context A. C\'ordoba and D. C\'ordoba \cite{cc} established the {\it a priori}  $L^p$ estimates: for $p\in[1,\infty]$ and $t\geq0$ 
\begin{equation}\label{eq-max}
\|\theta(t)\|_{L^p}\le\|\theta_0\|_{L^p}.
\end{equation}
 We remark that the proof in  the case $p=+\infty$  can be obtained from the following   representation of the fractional Laplacian $\DD^\beta$, 
 $$
 \DD^\beta f(x)=c_d\int_{\RR^d}\frac{f(x)-f(y)}{|x-y|^{d+\beta}}dy.
 $$
 Indeed, one can check  that if a continuous function reaches its  maximum at some  point $x_0$ then $\DD^\beta f(x_0)\geq0$ and hence we conclude as for the heat equation. 
 Our first main result si a generalization of the result  of \cite{cc} to \eqref{trdif}

\begin{Theo}\label{max-princ0} Let $\kappa\geq0, d\in\{2,3\},\beta\in]0,1],\alpha\geq0,\lambda\geq e^{\frac{3+2\alpha}{\beta}}$ and $  p\in[1,\infty]$. Then any smooth solution of \eqref{trdif} satisfies
$$
\|\theta(t)\|_{L^p} \le \|\theta_0\|_{L^p}.
$$

\end{Theo}
\begin{Rema}
The restriction on the parameter $\beta$ is technical and we believe that the above theorem remains true for $\beta\in]1,2[$ and $\alpha>0.$ 
\end{Rema}
Let us discuss the proof in the  special case of $v\equiv0.$ The equation \eqref{trdif} is reduced to the fractional heat equation

\centerline{$
\partial_{t}\theta+ \kappa\, \mathcal{L }\, \theta=0\quad\hbox{with}\quad\mathcal{L}:=\small{\frac{\DD^\beta}{\log^\alpha(\lambda+\DD)}}\cdot
$}
\

The solution is explicitly given by the convolution formula
$$
\theta(t,x)=K_t\star\theta_0(x)\quad\hbox{with}\quad \widehat{K_t}(\xi)=e^{-t \frac{|\xi|^\beta}{\log^\alpha(\lambda+|\xi|)}}.
$$ 
We will show that the family  $(K_t)_{t\geq0}$ is a convolution  semigroup of probabilities which means that  $\mathcal{L}$ is the  generator of a L\'evy semigroup. Consequently, this family    is a $C_0$-semigroup of contractions on $L^p$ for every $p\in[1,\infty[.$  The important  step in the proof   is to get  the positivity of  the kernel $K_t$. For this purpose we use Askey's criterion for  characteristic functions, see Theorem \ref{rad1}. We point out that the  restrictions on the dimension $d$  and the values of $\beta $ are due to the use of this criterion.
 Now to deal with the full transport-diffusion equation \eqref{trdif} we use some results from the theory  of $C_0-$semigroups of contractions. 
  
 The second estimate that we intend to establish is a generalized Bernstein inequality. Before stating the result we recall that for $q\in\NN$ the operator $\Delta_q$ is the frequency  localization around a ring of size $2^q$, see next section for more details.  Now our result reads as follows, \begin{Theo}\label{coer10}
Let $d\in\{1,2,3\},\beta\in]0,1],\alpha\geq0,\lambda\geq e^{\frac{3+2\alpha}{\beta}}$ and $  p\in]1,\infty[.$ Then we have for $q\in\NN$ and $f\in \mathcal{S}(\RR^d),$
\

\centerline{${
2^{q\beta}(q+1)^{-\alpha}\|\Delta_q f\|_{L^p}^p\le C{\displaystyle{\int}}_{\RR^d} \Big(\frac{\DD^{\beta}}{\log^{\alpha}(\lambda+\DD)}\Delta_qf\Big)\,|\Delta_q f|^{p-2}\Delta_qfdx.
}$}
\

where $C$ is a constant depending  on $p,\alpha,\beta$ and $\lambda$.
\end{Theo}
The proof relies on some tools from the theory of  L\'evy operators or more generally  sub-Markovians operators combined with some results from harmonic analysis. 
\begin{remas}
\begin{enumerate}
\item
When $\alpha=0$ then the above inequality is valid for all $\beta\in[0,2].$ The case $\beta=2$ was discussed in   \cite{dan,plan}. The   remaining case $\beta\in[0,2[$ was treated by Miao {\it et al.} in \cite{cmz} but only for $p\geq2.$
\item The proof for the case $p=2$  is an easy consequence of Plancherel identity and does not require any assumption on the parameters $ \alpha,\beta$ and $\lambda$.
\end{enumerate}

\end{remas}

The second part of this paper is concerned with   an application of Theorems \ref{max-princ0} and \ref{coer10} to the following  Boussinesq model with general dissipation
\begin{equation}
\label{bintro7}
\left\{ 
\begin{array}{ll} 
\partial_{t}v+v\cdot\nabla v+\nabla \pi=\theta e_{2},\,\, (t,x)\in\RR_+\times\RR^2  \\ 
\partial_{t}\theta+v\cdot\nabla\theta+ \kappa\, \mathcal{L}\,\theta=0\\
\textnormal{div}\,v=0\\
v_{| t=0}=v_{0}, \quad \theta_{| t=0}=\theta_{0}. 
\end{array} \right.
\end{equation}
Here, the velocity field $v$    is given by  $v=(v^1,v^2)$, the
  pressure $\pi$ and the temperature $\theta$ are scalar functions.
   The force term $\theta e_{2}$ in the velocity equation,  with $e_2$ the vector $(0,1)$, models
    the effect of the  gravity on  the fluid motion. The operator $\mathcal{L} $ 
    whose form may vary is used to  take into account 
    anomalous  diffusion  in the fluid motion.
  \
   From mathematical point of view, the question of global well-posedness for the  inviscid model, corresponding to $ \kappa = 0,$ is extremely hard to deal with. We point out that the classical theory of  symmetric hyperbolic quasi-linear systems can be applied for this system and thus  we can get the local well-posedness for smooth initial data. The significant quantity that one need to bound in order to get the global existence is the $L^\infty$-norm of the vorticity defined by  $\omega=
   \mbox{curl } v = \partial_{1} v^2 - \partial_{2} v^1.$ Now we observe from the first equation of \eqref{bintro7} that $\omega$ solves the equation
   \begin{equation}
   \label{vortintro}
   \partial_{t} \omega + v \cdot \nabla \omega =  \partial_{1} \theta.
   \end{equation}
   The main difficulty encountered for the global existence is due to the lack of strong dissipation in the temperature equation: we don't see how to estimate in a suitable way  the \mbox{quantity $\int_{0}^T \|\partial_{1} \theta \|_{L^\infty}$.}  
 However, the situation in the  viscous case, $\kappa>0$ and $ \mathcal{L}= -\Delta$, is well-understood and   the question of global existence is solved recently in a series of papers. In \cite{Cha}, Chae proved the global existence and uniqueness for initial data $(v_0,\theta_0)\in H^s\times H^s,$ with $s>2,$ see also  \cite{Hou}. This result was  improved  by the author and Keraani  in \cite{hk} 
   to initial data   $v_0\in B_{p,1}^{\frac{2}{p}+1}$ and $\theta_0\in B_{p,1}^{-1+\frac2p}\cap L^r, r>2$. The global existence of Yudovich solutions for this system was treated   in \cite{dp1}. We  also mention that in \cite{dp},  Danchin and Paicu  constructed global  strong   solutions for a dissipative term of the form $  \partial_{11}\theta$ instead of $\Delta\theta.$  In the papers \cite{H-Z,HKR} we try to understand the lower dissipation $\mathcal{L}=\DD^\alpha$ needed for global existence. In \cite{H-Z} the  author and Zerguine proved  the  global well-posedness when $\alpha\in]1,2[.$  The proof relies on   the fact that the dissipation is sufficiently strong to counterbalance the  possible amplification of  the vorticity due \mbox{to $\partial_1\theta.$} However the case $\alpha=1$ is not reached by the method and this value of $\alpha$ is  called critical   in the sense that  the dissipation and the amplification of the vorticity  due \mbox{to $\partial_{1} \theta $} have the same order.  \\
   
   In \cite{HKR} we prove that there is some hidden structure leading to global existence in the critical case. More precisely, we introduced the mixed quantity $\Gamma=\omega+\frac{\partial_1}{\DD}\theta$ which satisfies the equation
 $$
 \partial_t\Gamma+v\cdot\nabla\Gamma=-[\mathcal{R},v\cdot\nabla]\theta.
  \quad\hbox{with}\quad \mathcal{R}:={\frac{\partial_1}{\DD}}.$$
  As a matter of fact, the problem in the framework of Lebesgue spaces is reduced to the  estimate the commutator between the advection $v\cdot\nabla$ and Riesz transform $\mathcal{R}$ which is homogenous of degree zero. Since   Riesz transform is a Calder\'on-Zygmund operator then  using in a suitable way the  smoothing effects for $\theta$  we can get  a global estimate of $\|\omega(t)\|_{L^p}$. One can then use this information to control more strong norms of the vorticity like $\|\omega(t)\|_{L^\infty}$ or $\|\omega(t)\|_{B_{\infty,1}^0}$. For more discussions about the global well-posedness problem concerning  other classes of  Boussinesq  systems we refer to \cite{HKR1, MX}.
  
   Our goal here  is to relax the critical dissipation needed for global well-posedness by some logarithmic factor.    
    More precisely we will study the logarithmically critical Boussinesq model   
     
          \begin{equation}
     \label{Bouss}
\left\{ 
\begin{array}{ll} 
\partial_{t}v+v\cdot\nabla v+\nabla \pi=\theta e_{2}  \\ 
\partial_{t}\theta+v\cdot\nabla\theta +  \frac{\DD}{\log^\alpha(\lambda+\DD)} \theta = 0 \\
\textnormal{div}\,v=0\\
v_{| t=0}=v^{0}, \quad \theta_{| t=0}=\theta^{0}.  
\end{array} \right.
     \end{equation}

   Before stating our result we will need some new definitions.  First, we define the logarithmic Riesz transform $\mathcal{R}_\alpha$ by $\mathcal{R}_\alpha=\frac{\partial_1}{\DD}{\tiny}{\log^{\alpha}(\lambda+\DD)}.$ Second,  for a given $\alpha\in\RR$  we define the function spaces $\{\mathcal{X}_p\}_{1\le p\le\infty}$ by
   $$
   u\in \mathcal{X}_p\Leftrightarrow \|u\|_{\mathcal{X}_p}:=\|u\|_{B_{\infty,1}^0\cap L^p}+\|\mathcal{R}_\alpha u\|_{B_{\infty,1}^0\cap L^p}<\infty.
   $$
      
    Our result reads as follows (see section \ref{preliminaries} for the  definitions  and the basic properties of Besov spaces).
 \begin{Theo}\label{theo1}
   Let $\alpha\in[0,\frac12], \lambda\geq e^{4}$ and $p\in]2, \infty[$. Let $v_0\in {B}_{\infty, 1}^{1}\cap \dot{W}^{1,p}$ be a divergence free vector-field of $\RR^2$  and 
$\theta_0\in \mathcal{X}_p$. Then there exists a unique global solution $(v, \theta)$ to the system 
\eqref{Bouss} with 
\begin{equation*}
v\in L^\infty_{\textnormal{loc}}\big(\RR_+;{B}_{\infty, 1}^{1}\cap \dot{W}^{1,p}),
\qquad \theta\in L_{\textnormal{loc}}^{\infty}\big(\RR_+; \mathcal{X}_p\big)\cap \widetilde L_{\textnormal{loc}}^1(\RR_+; {B}_{p, \infty}^{1,-\alpha}).
\end{equation*}
\end{Theo}

The proof shares the same ideas as the case $\alpha=0$ treated in \cite{HKR} but with more technical difficulties. We define $\mathcal{R}_\alpha=\frac{\partial_1}{\DD}{\tiny}{\log^{\alpha}(\lambda+\DD)}$ and $\Gamma=\omega+\mathcal{R}_\alpha\theta.$ Then we get
$$
 \partial_t\Gamma+v\cdot\nabla\Gamma=-[\mathcal{R}_\alpha,v\cdot\nabla]\theta.
$$
To estimate the commutator in the framework of Lebesgue spaces we use the para-differential calculus combined with Theorems \ref{max-princ0}  and \ref{coer10}.
\begin{remas} 
\begin{enumerate}
\item We point out that for global well-posedness to the generalized Navier-Stokes system in dimension three, Tao proved in a recent paper \cite{Tao} that we can improve the dissipation $\DD^{\frac52}$   to $\frac{\DD^{\frac52}}{\log^{\frac12}(2+\DD)}\cdot$
\item The space $\mathcal{X}_ p$ is less regular than the space $B_{\infty,1}^\varepsilon\cap B_{p,1}^{\varepsilon}, \forall\varepsilon>0.$ More precisely, we will see in Proposition \ref{embed} that  $B_{\infty,1}^\varepsilon\cap B_{p,1}^{\varepsilon}\hookrightarrow \mathcal{X}_p.$
 \item If we take $\theta=0$ then the system \eqref{Bouss} is reduced to the  two-dimensional  Euler system.  It is well known that this system is globally well-posed in  $H^s$ for  $s>2$. The main tool for global existence is the  BKM criterion  \cite{bkm} ensuring that the development of finite-time singularities for Kato's solutions is related to the blowup  of the $L^\infty$ norm of the vorticity near the maximal time existence.
In \cite{vis} Vishik    extended the  global existence  of strong solutions   to initial data belonging to Besov  spaces  $B_{p,1}^{1+2/p}.$  
Notice that  these spaces   have the same scale as  Lipschitz functions   and in this sense they are called  critical and  it is not at all clear whether BKM criterion can be used in this context.    
\item
Since 
$
 B_{r,1}^{1+2/r}\hookrightarrow {B}_{\infty, 1}^{1}\cap \dot{W}^{1,p}$ for all $ r\in [1,+\infty[$ and $p>\max\{r,2\}$, then 
  the space of initial velocity in our theorem  contains all the critical spaces $B_{p,1}^{1+2/p}$ except the biggest  one, that is  $B_{\infty,1}^{1}$. 
 For the limiting case  we have been  able to prove the global existence only  up to  the extra assumption   $\nabla v_0\in L^p$ for some $p\in]2, \infty[$.  The reason behind this extra assumption is the fact that to  obtain a global  $L^\infty$ bound for the vorticity we need before to establish an $L^p$ estimate for some $p\in]2,\infty[$ and it is not clear how to get rid of  this condition.
\item
 Since  $\nabla v\in L^{1}_{\rm loc}(\RR_+;L^\infty)$   then we can propagate all the higher regularities: critical (for example  $v_0\in B_{p,1}^{1+2/p}$ with $p<\infty$) and sub-critical (for example  $v_0\in H^{s}$, \mbox{with $s>2$).} \end{enumerate}
\end{remas}

\section{Notations and preliminaries}
    
    \label{preliminaries}
       \subsection{Notations}Throughout this paper we will use the following notations.
       
$\bullet$ For any positive  $A$ and $B$  the notation  $A\lesssim B$ means that there exist a positive  harmless constant $C$ such that $A\le CB$. 

$\bullet$ For any tempered distribution $u$  both $ \hat u$  and $\mathcal F u$ denote the Fourier transform of $u$.

$\bullet$ Pour every $p\in [1,\infty]$, $\|\cdot\|_{L^p}$ denotes the norm in the Lebesgue space $L^p$.

$\bullet$ The norm  in the mixed space time Lebesgue space $L^p([0,T],L^r(\mathbb R^d)$ is denoted by  $\|\cdot\|_{L^p_TL^r}$ (with the obvious generalization to  $\|\cdot\|_{L^p_T\mathcal X} $ for any normed space $\mathcal X$).

$\bullet$ For any pair of operators $P$ and $Q$ on some Banach space $\mathcal{X}$, the commutator $[P,Q]$ is given by $PQ-QP$.

$\bullet$ For $p\in[1,\infty]$, we denote by $\dot{W}^{1,p}$ the space of distributions $u$ such that $\nabla u\in L^p.$

    \subsection{Functional spaces}

  Let us introduce the so-called   Littlewood-Paley decomposition and the corresponding  cut-off operators. 
There exists two radial positive  functions  $\chi\in \mathcal{D}(\RR^d)$ and  $\varphi\in\mathcal{D}(\RR^d\backslash{\{0\}})$ such that
\begin{itemize}
\item[\textnormal{i)}]
$\displaystyle{\chi(\xi)+\sum_{q\geq0}\varphi(2^{-q}\xi)=1}$;$\quad \displaystyle{\forall\,\,q\geq1,\, \textnormal{supp }\chi\cap \textnormal{supp }\varphi(2^{-q})=\varnothing}$\item[\textnormal{ii)}]
 $ \textnormal{supp }\varphi(2^{-j}\cdot)\cap
\textnormal{supp }\varphi(2^{-k}\cdot)=\varnothing,$ if  $|j-k|\geq 2$.
\end{itemize}

For every $v\in{\mathcal S}'(\RR^d)$ we set 
  $$
\Delta_{-1}v=\chi(\hbox{D})v~;\, \forall
 q\in\NN,\;\Delta_qv=\varphi(2^{-q}\hbox{D})v\quad\hbox{ and  }\;
 S_q=\sum_{ j=-1}^{ q-1}\Delta_{j}.$$
 The homogeneous operators are defined by
 $$
 \dot{\Delta}_{q}v=\varphi(2^{-q}\hbox{D})v,\quad \dot S_{q}v=\sum_{j\leq q-1}\dot\Delta_{j}v,\quad\forall q\in\ZZ.
 $$
From   \cite{b}  we split the product 
  $uv$ into three parts: $$
uv=T_u v+T_v u+R(u,v),
$$
with
$$T_u v=\sum_{q}S_{q-1}u\Delta_q v,\quad  R(u,v)=\sum_{q}\Delta_qu\tilde\Delta_{q}v  \quad\hbox{and}\quad \tilde\Delta_{q}=\sum_{i=-1}^1\Delta_{q+i}.
$$

\

 For
 $(p,r)\in[1,+\infty]^2$ and $s\in\RR$ we define  the inhomogeneous Besov \mbox{space $B_{p,r}^s$} as 
the set of tempered distributions $u$ such that
$$\|u\|_{B_{p,r}^s}:=\Big( 2^{qs}
\|\Delta_q u\|_{L^{p}}\Big)_{\ell ^{r}}<+\infty.
$$
The homogeneous Besov space $\dot B_{p,r}^s$ is defined as the set of  $u\in\mathcal{S}'(\RR^d)$ up to polynomials such that
$$
\|u\|_{\dot B_{p,r}^s}:=\Big( 2^{qs}
\|\dot\Delta_q u\|_{L^{p}}\Big)_{\ell ^{r}(\ZZ)}<+\infty.
$$
For $s,s'\in\RR$ and $p,r\in[1,\infty]$ we define the generalized Besov space $ B_{p,r}^{s,s'}$ as the set of tempered distributions $u$ such that
$$
\|u\|_{B_{p,r}^{s,s'}}:=\Big(2^{qs}(|q|+1)^{s'}\|\Delta_q u\|_{L^p}\Big)_{\ell^r}<\infty
$$
{Let $T>0$} \mbox{and $\rho\geq1,$} we denote by $L^\rho_{T}B_{p,r}^{s,s'}$ the space of distributions $u$ such that 
$$
\|u\|_{L^\rho_{T}B_{p,r}^{s,s'}}:= \Big\|\Big( 2^{qs}(|q|+1)^{s'}
\|\Delta_q u\|_{L^p}\Big)_{\ell ^{r}}\Big\|_{L^\rho_{T}}<+\infty.$$
We say that 
$u$ belongs to the space
 $\widetilde L^\rho_{T}{B_{p,r}^{s,s'}}$ if
 $$
 \|u\|_{ \widetilde L^\rho_{T}{B_{p,r}^s}}:= \Big( 2^{qs}(|q|+1)^{s'}
\|\Delta_q u\|_{L^\rho_{T}L^p}\Big)_{\ell ^{r}}<+\infty.$$
By  a direct application  of 
the Minkowski inequality, we have the following links between these spaces. 
Let $ \varepsilon>0,$ then 
$$
L^\rho_{T}B_{p,r}^s\hookrightarrow\widetilde L^\rho_{T}{B_{p,r}^s}\hookrightarrow{L^\rho_{T}}{B_{p,r}^{s-\varepsilon}}
,\,\textnormal{if}\quad  r\geq \rho,$$
$$
{L^\rho_{T}}{B_{p,r}^{s+\varepsilon}}\hookrightarrow\widetilde L^\rho_{T}{B_{p,r}^s}\hookrightarrow L^\rho_{T}B_{p,r}^s,\, \textnormal{if}\quad 
\rho\geq r.
$$

 We will  make continuous use of Bernstein inequalities (see  \cite{che1} for instance).
\begin{lem}\label{lb}\;
 There exists a constant $C$ such that for $q,k\in\NN,$ $1\leq a\leq b$ and for  $f\in L^a(\RR^d)$, 
\begin{eqnarray*}
\sup_{|\alpha|=k}\|\partial ^{\alpha}S_{q}f\|_{L^b}&\leq& C^k\,2^{q(k+d(\frac{1}{a}-\frac{1}{b}))}\|S_{q}f\|_{L^a},\\
\ C^{-k}2^
{qk}\|{\Delta}_{q}f\|_{L^a}&\leq&\sup_{|\alpha|=k}\|\partial ^{\alpha}{\Delta}_{q}f\|_{L^a}\leq C^k2^{qk}\|{\Delta}_{q}f\|_{L^a}.
\end{eqnarray*}

\end{lem}

\section{Maximum principle}
Our task is to establish some useful estimates for the following equation generalizing \eqref{trdif} 
\begin{equation}
      \label{transport-d} 
      \left\{ \begin{array}{ll} 
\partial_t\theta+v\cdot\nabla\theta+\frac{\DD^\beta}{\log^\alpha(\lambda+\DD)}\theta=f\\
\textnormal{div}\,  v=0\\
 \theta_{| t=0}=\theta_{0},
\end{array} \right.
\end{equation}
Two special  problems will be studied: the first one deals with  $L^p$ estimates that give in particular Theorem \ref{max-princ0}. However the second one  consists in establishing some logarithmic estimates in Besov spaces with index regularity zero.  
 The first main result of this section generalizes Theorem \ref{max-princ0}.
\begin{Theo}\label{max-princ} Let $p\in[1,\infty],\beta\in]0,1],\alpha\geq0$ and $\lambda\geq e^{\frac{3+2\alpha}{\beta}}$. Then any smooth solution of \eqref{transport-d} satisfies
$$
\|\theta(t)\|_{L^p}\le\|\theta_0\|_{L^p}+\int_0^t\|f(\tau)\|_{L^p}d\tau.
$$

\end{Theo}

The proof will be done in two steps. The first one is to valid the result for the free fractional heat equation. More precisely we will establish that the semigroup $e^{t\mathcal{L}}$, with $\mathcal{L}:=\frac{\DD^\beta}{\log^\alpha(\lambda+\DD)}$, is a contraction in Lebesgue spaces $L^p$, for every $p\in[1,\infty[$ of course under suitable conditions on the parameters $\alpha,\beta,\lambda$.
This problem is reduced to show that 
$\|K_t\|_{L^1}\le1$. This  is equivalent to $K_t\in L^1$ and $K_t\geq0.$ As we will see, to get the integrability of the kernel  we do not need any restriction on the value of our parameters. Nevertheless,  the positivity of $K_t$ requires some restrictions which are detailed in \mbox{Theorem \ref{max-princ}.} The second step is to establish the $L^p$ estimate for the system  \eqref{transport-d} and for this purpose we use some results about L\'evy operators or more generally sub-Markovians operators. 
\subsection{Definite positive functions}
As we will see there is a strong connection between the  positivity of the kernel $K_t$ introduced above  and the notion of definite positive functions. We will first gather some well-known properties about definite positive functions and recall some useful criteria for characteristic functions. Second and as an application we will show that the kernel $K_t$ is positive under suitable conditions on the involved parameters.
\begin{defi}
Let $f:\RR^d\to\CC$ be a complex-valued function. We say that $f$ is { definite positive} if only if for every integer $n\in\NN^*$ and every set of points  $\{x_j,j=1,...,n\}$ of $\RR^d$ the matrix $(f(x_j-x_k))_{1\le j,k\le n}$ is  positive Hermitian, that is, for every complex numbers $\xi_1,...,\xi_n$ we have
$$
\sum_{j,k=1}^nf(x_j-x_k)\xi_j\bar\xi_k\geq 0.
$$
\end{defi}
We will give some  results about  definite positive functions.
\begin{enumerate}
\item
From the definition, every  definite positive function $f$ satisfies
$$
f(0)\geq 0,\quad f(-x)=\overline{f(x)},\quad |f(x)|\le f(0).
$$
\item The continuity of a definite positive function $f$ at zero gives the continuity everywhere. More precisely we have
$$
|f(x)-f(y)|\le 2f(0)\big(f(0)-f(x-y)\big).
$$
\item The sum of two  definite positive functions is also definite positive and according to  Shur's lemma  the product of two  definite positive functions  is also  definite positive and therefore  the class of definite positive functions is a convex cone closed under multiplication. 

\item Let $\mu$ be a finite positive measure then its Fourier-Stielitjes transform is given by
$$
\widehat\mu(\xi)=\int_{\RR^d}e^{-ix\cdot\xi}d\mu(x).
$$
It is easy to see that $\widehat{\mu}$ is a   definite positive function.  Indeed
\begin{eqnarray*}
\sum_{j,k=1}^n\widehat\mu(x_j-x_k)\xi_j\bar \xi_k&=&\int_{\RR^d}\Big(\sum_{j,k=1}^ne^{-ix\cdot x_j}\xi_je^{ix\cdot x_k}\bar \xi_k\Big)d\mu(x)\\
&=&\int_{\RR^d}\Big|\sum_{j=1}^ne^{-ix\cdot x_j}\xi_j\Big|^2d\mu(x)\\
&\geq&0.
\end{eqnarray*}
\end{enumerate}

The converse  of the last point $(4)$ is stated by the following result due to Bochner, see for instance Theorem 19 in  \cite{Bochner}. 
\begin{Theo}[Bochner's theorem] Let $f:\RR^d\to\CC$ be a continuous definite positive  function, then $f$ is the Fourier transform of a finite positive  Borel measure.

\end{Theo}
Hereafter we will focus on a the class of radial definite positive functions. First we say that $f:\RR^d\to\RR$ is radial if $f(x)=F(|x|)$ with 
$F:[0,+\infty[\to\CC$. 
There are some criteria  for radial functions to be definite positive. For example in dimension one the celebrated criterion of  P\'olya \cite{Polya} states that  if $F:[0,+\infty[\to\RR$  is continuous and  convex  with $F(0)=1$ and $\lim_{r\to+\infty}F(r)=0$ then $f(x)=F(|x|) $ is  definite positive. This criterion was extended in higher dimensions by numerous authors \cite{Ask,Gneit,Trig}. We will restrict ourselves to the following one due to Askey \cite{Ask}.
\begin{Theo}[Askey] \label{rad1} Let $d\in\NN,\,\,F:[0,+\infty\to\RR$ be a continuous  function such that
\begin{enumerate}
\item 
$F(0)=1,$ 
\item
the function $r\mapsto (-1)^dF^{(d)}(r)$ exists and is convex on $]0,+\infty[,$
\item $lim_{r\to+\infty}F(r)=\lim_{r\to+\infty} F^{(d)}(r)=0.$ 

\end{enumerate} 
Then for every $k\in\{1,2,..,2d+1\} $ the function $x\mapsto F(|x|)$ is the Fourier transform of a probability measure on $\RR^k$.

\end{Theo}
\begin{Rema}
As an  application of Askey's theorem we have that $x\mapsto e^{-t|x|^\beta}$ is definite positive for all $t>0,$ $\beta\in]0,1]$ and $d\in\NN$. Indeed, the function $F(r)=e^{-tr^\beta}$ is completely monotone, that is,  
$(-1)^kF^{(k)}(r)\geq0,\,\forall \,r>0,\, k\in\NN.$ Although the case $\beta\in]1,2]$ can not be reached by this criterion the result is still true.

\end{Rema} We will now see that the perturbation of the above function by a logarithmic damping is also definite positive.  More precisely, we have
\begin{prop}\label{rad2}
Let $\alpha,t\in[0,+\infty[\times]0,+\infty[, \beta\in]0,1],\,\lambda\geq e^{\frac{3+2\alpha}{\beta}}$ and define $f:\RR^d\to\RR$  by
$$
f(x)=e^{-t\frac{|x|^\beta}{\log^{\alpha}(\lambda+|x|)}}.
$$
Then $f$ is a definite positive function  for $d\in\{1,2,3\}.$
\end{prop}
\begin{remas}
\begin{enumerate}
\item It is possible that the above result remains true for higher dimension $d\geq4$ but we avoid to deal with this more computational case. We think also that the radial function associated to $f$ is completely monotone.
\item The lower bound of $\lambda$ is not optimal by our method. In fact we can obtain more precise bound but this seems to be irrelevant.
\end{enumerate}
\end{remas}
\begin{proof}
We  write  $f(x)=F(|x|)$ with 
$$F(r)=e^{-t\phi(r)}\quad\hbox{and}\quad \phi(r)=\frac{r^\beta}{\log^{\alpha}(\lambda+r)}\cdot
$$ 
 The function $F$ is smooth on $]0,\infty[$ and assumptions $(1)$ and $(3)$ of Theorem \ref{rad1} are satsified. It follows that  the function $f$ is definite positive  for $d\in\{1,2,3\}$  if 
$$F^{(3)}(r)\le0.
$$ Easy computations give for $r>0$,
$$
F^{(3)}(r)=\Big[-t\,\phi^{(3)}(r)+3t^2\,\phi^{\prime}(r)\phi^{(2)}(r)-t^3\,(\phi^{\prime}(r))^3\Big]F(r).
$$
We will prove that
$$
\phi^{\prime}(r)\geq0,\, \phi^{(2)}(r)\le0\quad\hbox{and}\quad \phi^{(3)}(r)\geq0.
$$
This is sufficient to get  $F^{(3)}(r)\le0,\forall\, r>0.$ The first derivative of $\phi$ is given by
\begin{eqnarray*}
\phi^{\prime}(r)&=&\frac{\beta\, r^{\beta-1}}{\log^{\alpha}(\lambda+r)}-\frac{\alpha\, r^{\beta}}{(\lambda+r)\log^{\alpha+1}(\lambda+r)}\\
&=&\frac{ r^{\beta-1}}{(\lambda+r)\log^{\alpha+1}(\lambda+r)}\Big(\beta\lambda\log(\lambda+r)+r\big(\beta\log(\lambda+r)-\alpha\big)   \Big).
\end{eqnarray*}
We see that if $\lambda$ satisfies
\begin{equation}\label{c1}
\lambda\geq e^{\frac\alpha\beta}
\end{equation}
then $\phi^{\prime}(r)\geq0$. For the second derivative of $\phi$ we obtain
\begin{eqnarray*}
\phi^{(2)}(r)&=&-
\frac{\beta(1-\beta)r^{\beta-2}}{\log^{\alpha}(\lambda+r)}-\frac{2\alpha\beta r^{\beta-1}}{(\lambda+r)\log^{1+\alpha}(\lambda+r)}\\
&+&\frac{\alpha r^\beta}{(\lambda+r)^2\log^{\alpha+1}(\lambda+r)}+\frac{\alpha(\alpha+1)r^\beta}{(\lambda+r)^2\log^{\alpha+2}(\lambda+r)}\\
&=&\frac{r^{\beta-2}}{\log^{\alpha}(\lambda+r)}\Bigg[ -\beta(1-\beta)-\frac{2\alpha\beta r}{(\lambda+r)\log(\lambda+r)}\\
&+&\frac{\alpha r^2}{(\lambda+r)^2\log(\lambda+r)}+\frac{\alpha(\alpha+1) r^2}{(\lambda+r)^2\log^2(\lambda+r)}\Bigg]. 
\end{eqnarray*}
Since $\frac{ r^2}{(\lambda+r)^2}\le\frac{ r}{\lambda+r}\le1,$ then 
\begin{eqnarray*}
\phi^{(2)}(r)&\le&\frac{r^{\beta-2}}{\log^{\alpha}(\lambda+r)}\Bigg[(1-\beta)\Big(-\beta+\frac{2\alpha}{\log(\lambda+r)}\Big)-\frac{\alpha r}{(\lambda+r)\log(\lambda+r)}\Big(1-\frac{\alpha+1}{\log(\lambda+r)}\Big)\Bigg]\\
&\le&\frac{r^{\beta-2}}{\log^{\alpha}(\lambda+r)}\Bigg[(1-\beta)\Big(-\beta+\frac{2\alpha}{\log\lambda}\Big)-\frac{\alpha r}{(\lambda+r)\log(\lambda+r)}\Big(1-\frac{\alpha+1}{\log\lambda}\Big)\Bigg].
\end{eqnarray*}
Now we choose $\lambda$ such that 
$$
-\beta+\frac{2\alpha}{\log\lambda}\leq0\quad\hbox{and}\quad 1-\frac{\alpha+1}{\log\lambda}\geq0
$$
which is true despite $\lambda$ satisfies
\begin{equation}\label{cond1}
\max(e^{\frac{2\alpha}{\beta}},e^{\alpha+1})\le\lambda.
\end{equation}
Under this assumption we get 
$$\phi^{(2)}(r)\le0,\forall r>0.
$$
Similarly we have
\begin{eqnarray*}
\phi^{(3)}(r)&=&\alpha(\alpha+1)r^{\beta-1}\frac{\log^{-\alpha-3}(\lambda+r)}{(\lambda+r)^2}\Big[3\lambda\beta\log(\lambda+r)+r\big(3\beta\log(\lambda+r)-(2+\alpha)\big)\Big]\\&+&
\alpha r^{\beta-2}\frac{\log^{-2-\alpha}(\lambda+r)}{(\lambda+r)^3}\Big[r^2\Big(-3(1+\alpha)+(-3\beta^2+6\beta-2)\log(\lambda+r)\Big)\\
&+&\log(\lambda+r)\Big(\lambda\beta(9-6\beta)r+3\lambda^2\beta(1-\beta)  \Big)  \Big]\\&+&
(2-\beta)(1-\beta)\beta r^{\beta-3}\log^{-\alpha}(\lambda+r)\\&=&
I_1+I_2+I_3+I_4.
\end{eqnarray*}
it is easy to see that $I_3$ and $I_4$ are nonnegative. On the other hand we have
\begin{eqnarray*}
I_1+I_2&=&3\lambda\beta\alpha(\alpha+1)r^{\beta-1}\frac{\log^{-\alpha-2}(\lambda+r)}{(\lambda+r)^2}\\&+&
\alpha r^{\beta}\frac{\log^{-2-\alpha}(\lambda+r)}{(\lambda+r)^3}\Bigg[-3(1+\alpha)+(\alpha+1)(\lambda+r)\Big(3\beta-\frac{2+\alpha}{\log(\lambda+r)}\Big)\\&+&(-3\beta^2+6\beta-2)\log(\lambda+r)\Bigg]
\end{eqnarray*}
Since $-3\beta^2+6\beta-2\geq-2,$  for $\beta\in[0,1]$  and $-\frac{\log x}{x}\geq -\frac{\log\lambda}{\lambda}, \forall x\geq \lambda\geq e,$ then
\begin{eqnarray*}
I_1+I_2&\geq&
\alpha(\alpha+1) r^{\beta}\frac{\log^{-2-\alpha}(\lambda+r)}{(\lambda+r)^3}\Bigg[-3+(\lambda+r)\Big(3\beta-\frac{2+\alpha}{\log\lambda}-2\frac{\log(\lambda+r)}{(\alpha+1)(\lambda+r)}\Big)\Bigg]\\
&\geq&\alpha(\alpha+1) r^{\beta}\frac{\log^{-2-\alpha}(\lambda+r)}{(\lambda+r)^2}\Bigg[3\beta-\frac{3}{\lambda}-\frac{2+\alpha}{\log\lambda}-\frac{2\log\lambda}{(\alpha+1)\lambda}\Big)\Bigg].
\end{eqnarray*}
We can check that 
$$
\log\lambda\le\lambda\quad\hbox{and}\quad \log^2\lambda\le\lambda, \,\forall\lambda\geq e.
$$
Thus
\begin{eqnarray*}
I_1+I_2&\geq&\alpha(\alpha+1) r^{\beta}\frac{\log^{-2-\alpha}(\lambda+r)}{(\lambda+r)^2}\Big[3\beta-\frac{1}{\log\lambda}\big(5+\alpha+\frac{2}{\alpha+1}  \big)\Big]\\
&\geq&\alpha(\alpha+1) r^{\beta}\frac{\log^{-2-\alpha}(\lambda+r)}{(\lambda+r)^2}\Big[3\beta-\frac{7+\alpha}{\log\lambda}\Big].
\end{eqnarray*}
We choose $\lambda$ such that
\begin{equation*}
3\beta-\frac{7+\alpha}{\log\lambda}\geq 0.
\end{equation*}
 It follows that  $I_1+I_2$ is nonnegative if
\begin{equation}\label{c3}
\lambda\geq e^{\frac{7+\alpha}{3\beta}}.
\end{equation}
Remark that the assumptions \eqref{c1}, \eqref{cond1} and \eqref{c3} are satisfies under the condition 
$$
\lambda\geq e^{\frac{3+2\alpha}{{\beta}}}.
$$
Finally, we get: $\quad \forall \alpha\geq0,\, \beta\in]0,1],\,\lambda\geq e^{\frac{3+2\alpha}{{\beta}}},$
$$
 \forall r>0,\quad\phi^{(3)}(r)\geq0.
$$
This achieves the proof.
\end{proof}
More  precise informations  about the kernel $K_t$ will be listed in the following lemma. \begin{lem}\label{kernel}
 Let $\lambda\geq 2$ and denote by  $K_t$ the element of $\mathcal{S}'(\RR^d)$  such that 
 $$\widehat{K_t}(\xi)=e^{-t\frac{|\xi|^\beta}{\log^{\alpha}(\lambda+|\xi|)}}.
 $$ Then we have the following properties. 
\begin{enumerate}
\item For $(t,\alpha,\beta)\in]0,\infty[\times\RR\times]0,\infty[$ the function $K_t$ belongs to $\in L^1\cap C_0$.
\item For  $d\in\{1,2,3\},(t,\alpha,\beta)\in]0,+\infty[\times[0,\infty[\times]0,1]$ and $\lambda\geq e^{\frac{3+2\alpha}{\beta}},$ we have 
$$
\quad K_t(x)\geq 0,\,\forall x\in \RR_+\quad\hbox{and}\quad \|K_t\|_{L^1}=1.
$$
\end{enumerate}
\end{lem}

\begin{proof}
{\bf$(1)$} 
 By definition we have
$$
K_t(x)=(2\pi)^{-d}\int_{\RR^d}e^{-t\frac{|\xi|^\beta}{\log^{\alpha}(\lambda+|\xi|)}} e^{i\,x\cdot\xi}d\xi.
$$
Let $\mu\geq0,$ then integrating by parts we get
$$
|x|^{\mu}x_j^dK_t(x)=(-2i\pi)^{-d}\int_{\RR^d}\partial_{\xi_j}^d\big(e^{-t\frac{|\xi|^\beta}{\log^{\alpha}(\lambda+|\xi|)}}\big) |x|^{\mu}e^{i\,x\cdot\xi}d\xi.
$$
On the other hand we have
$$
|x|^{\mu}e^{i\,x\cdot\xi}=\DD^{\mu}e^{i\,x\cdot\xi},
$$
here $\DD$ is a fractional derivative on the variable $\xi$. Thus we get
$$
|x|^{\mu}x_j^dK_t(x)=(-2i\pi)^{-d}\int_{\RR^d}\DD^\mu\partial_{\xi_j}^d\big(e^{-t\frac{|\xi|^\beta}{\log^{\alpha}(\lambda+|\xi|)}}\big)e^{i\,x\cdot\xi}d\xi.
$$
Now we use the following representation for $\DD^\mu$ when $\mu\in]0,2]$
$$
\DD^\mu f(x)=C_{\mu,d}\int_{\RR^{d}}\frac{f(x)-f(x-y)}{|y|^{d+\mu}}dy.
$$
It follows that
$$
|x|^{\mu}|x_j^dK_t(x)|\le C_{\mu,d}\int_{\RR^{2d}}\frac{|\mathcal{K}_j(\xi)-\mathcal{K}_j(\xi-y)|}{|y|^{d+\mu}}dyd\xi
$$
with 
$$
\mathcal{K}_j(\xi):=\partial_{\xi_j}^d\big(e^{-t\frac{|\xi|^\beta}{\log^{\alpha}(\lambda+|\xi|)}}\big).
$$
Now we decompose the integral into two parts
\begin{eqnarray*}
\int_{\RR^{2d}}\frac{|\mathcal{K}_j(\xi)-\mathcal{K}_j(\xi-y)|}{|y|^{d+\mu}}dyd\xi&=&\int_{|y|\geq\frac{|\xi|}{2}}\frac{|\mathcal{K}_j(\xi)-\mathcal{K}_j(\xi-y)|}{|y|^{d+\mu}}dyd\xi\\
&+&\int_{|y|\leq\frac{|\xi|}{2}}\frac{|\mathcal{K}_j(\xi)-\mathcal{K}_j(\xi-y)|}{|y|^{d+\mu}}dyd\xi\\
&=&I_1+I_2.
\end{eqnarray*}
To estimate the first term we use the  following estimate that can be obtained by straightforward computations
\begin{eqnarray*}
|\mathcal{K}_j(\xi)|&\le& C_{t,\alpha,\beta}\frac{|\xi|^{\beta-d}}{\log^{\alpha}(\lambda+|\xi|)}e^{-t\frac{|\xi|^\beta}{\log^{\alpha}(\lambda+|\xi|)}}\\
&\le&C_{t,\alpha,\beta}{|\xi|^{\beta-d}}e^{-\frac12t\frac{|\xi|^\beta}{\log^{\alpha}(\lambda+|\xi|)}}.
\end{eqnarray*}
Hence we get under the assumption $\mu\in]0,\beta[,$
\begin{eqnarray*}
I_1&\le& C_{t,\alpha,\beta}\int_{|\xi|\le 2|y|}\frac{1}{|y|^{d+\mu}}\Big({|\xi|^{\beta-d}}e^{-\frac12t\frac{|\xi|^\beta}{\log^{\alpha}(\lambda+|\xi|)}}+{|\xi-y|^{\beta-d}}e^{-\frac12t\frac{|\xi-y|^\beta}{\log^{\alpha}(\lambda+|\xi-y|)}}\Big)d\xi dy\\
&\le&C_{t,\alpha,\beta}\int_{|\xi|\le 3|y|}\frac{1}{|y|^{d+\mu}}{|\xi|^{\beta-d}}e^{-\frac12 t\frac{|\xi|^\beta}{\log^{\alpha}(\lambda+|\xi|)}}d\xi dy\\
&\le&C_{t,\alpha,\beta}\int_{\RR^d}\frac{1}{|\xi|^{d+\mu-\beta}}e^{-\frac12 t\frac{|\xi|^\beta}{\log^{\alpha}(\lambda+|\xi|)}}d\xi\\
&\le& C_{t,\alpha,\beta}.
\end{eqnarray*}
To estimate the second term we use the mean-value Theorem
$$
{|\mathcal{K}_j(\xi)-\mathcal{K}_j(\xi-y)|}\le|y|\sup_{\eta\in[\xi-y,\xi]}|\nabla \mathcal{K}_j(\eta)|.
$$
On the other hand we have
\begin{eqnarray*}
|\nabla \mathcal{K}_j(\eta)|&\le& C_{t,\alpha,\beta}|\eta|^{\beta-d-1}e^{-\frac12t\frac{|\eta|^\beta}{\log^{\alpha}(\lambda+|\eta|)}}.
\end{eqnarray*}
Now since $|y|\le\frac12|\xi|$ then for $\eta\in[\xi-y,\xi]$ we have 
$$
\frac12|\xi|\le|\eta|\le\frac52|\xi|.
$$
This yields 
$$
{|\mathcal{K}_j(\xi)-\mathcal{K}_j(\xi-y)|}\le C_t|y||\xi|^{\beta-d-1}e^{-Ct|\xi|^{\frac\beta2}}.
$$
Therefore we find for $\mu\in]0,\beta[\cap]0,1[$,
\begin{eqnarray*}
I_2&\le& C_{t,\alpha,\beta}\int_{|y|\le\frac12|\xi|}\frac{1}{|y|^{d+\mu-1}}|\xi|^{\beta-d-1}e^{-Ct|\xi|^{\frac\beta2}}dyd\xi\\
&\le&C_{t,\alpha,\beta}\int_{\RR^2}\frac{1}{|\xi|^{d+\mu-\beta}}e^{-Ct|\xi|^{\frac\beta2}}d\xi\\
&\le&C_{t,\alpha,\beta}.
\end{eqnarray*}
Finally we get
$$
\hbox{for}\quad j=1,..,d,\quad |x|^\mu|x_j|^d|K_t(x)|\le C_{t,\alpha,\beta}.
$$
Since $K_t\in C_0$ then 
$$
\big(1+|x|^{d+\mu}\big)|K_t(x)|\le C_t.
$$
This proves that $K_t\in L^1(\RR^d)$. 

{\bf $(2)$} Using Theorem \ref{rad1} and Proposition \ref{rad2} we get $K_t\geq0$. Since $K_t\in L^1$ then 
$$
\|K_t\|_{L^1}=\widehat{K_t}(0)=1.
$$
\end{proof}
 Now we define 
 the propagator $e^{-t\frac{\DD^\beta}{\log^\alpha(\lambda+\DD)}}$ by  convolution 
 $$
 e^{-t\frac{\DD^\beta}{\log^\alpha(\lambda+\DD)}}f=K_t\star f.
 $$
 We have the following result.
\begin{cor}\label{cor1}
Let $\alpha\geq0,\beta\in]0,1],\,\lambda\geq e^{\frac{3+2\alpha}{\beta}}$ and $p\in[1,\infty]$. Then 
$$
\|e^{-t\frac{\DD^\beta}{\log^\alpha(\lambda+\DD)}} f\|_{L^p}\le\|f\|_{L^p},\quad\forall\,\, f\in L^p.
$$

\end{cor}
\begin{proof}
From the classical convolution inequalities combined with Lemma \ref{kernel} we get
\begin{eqnarray*}
\|e^{-t\frac{\DD^\beta}{\log^\alpha(\lambda+\DD)}} f\|_{L^p}&\le&\|K_t\|_{L^1}\|f\|_{L^p}\\
&\le&\|f\|_{L^p}.
\end{eqnarray*}

\end{proof}

\subsection{Structure of the semigroup $(e^{-t\frac{\DD^\beta}{\log^\alpha(\lambda+\DD)}})_{t\geq0}$}

We will first recall the notions of $C_0-$semigroup and  sub-Markovian generators. First we introduce the notion of strongly continuous semigroup.
\begin{defi}\label{def-sem}
Let $X$ be a Banach space and $(T_t)_{t\geq0}$ be  a family of bounded operators from $X$ into $X$. This family  is called a strongly continuous semigroup on $X$ or a $C_0$-semigroup if
\begin{enumerate}
\item $T_0=\hbox{Id},$
\item for every $t,s\geq0,\, T_{t+s}=T_tT_s,$
\item for every $x\in X, \lim_{t\to 0^+}\|T_t x-x\|=0.$ 
\end{enumerate}
If in addition the semigroup satisfies the estimate 
$$
\|T_t\|_{\mathcal{L}(X)}\leq 1,
$$
then it is called a $C_0$-semigroup of contractions.
\end{defi}
For a given $C_0$-semigroup of contractions  $(T_t)_{t\geq0}$ we define its domain $\mathcal{D}(A)$ by
$$
\mathcal{D}(A):=\Big\{f\in X;\lim_{t\to0^+}\frac{T_tf-f}{t}\,\,\textnormal{exists in}\,\, X \Big\},
$$ 

$$
Af=\lim_{t\to0^+}\frac{T_tf-f}{t},\quad f\in \mathcal{D}(A).
$$
It is well-known that the operator $A$ is densely defined: its  domain $\mathcal{D}(A)$ is dense \mbox{in $X.$}

We introduce now the special case of sub-Markovian semigroups.
\begin{defi}\label{defis2}
Let $X=L^p(\RR^d),$ with $p\in[1,\infty[ $ and $d\in\NN^*.$ Let $(T_t)_{t\geq0}$ be a $C_0$-semigroup of contractions on $X.$ It is said a sub-Markovian semigroup  if
\begin{enumerate}
\item If $f\in X, f(x)\geq 0,$ a .e.  then for every $t\geq0$,  $T_t f(x)\geq 0,$ a.e.,
\item  If $f\in X, |f|\leq1$ then for every $t\geq0,\,|T_tf|\le1$.
\end{enumerate}
\end{defi}
Denote by $L^p_+:=\Big\{f\in L^p; f(x)\geq0, a.e\Big\}$. Then we have the following classical result.
\begin{Theo}[Beurling-Deny theorem]\label{deurl}
Let $A$ be a nonnegative self-adjoint operator \mbox{of $L^2$.} Then we have the equivalence between
\begin{enumerate}
\item $\forall t>0,\, f\in L^2_+\Rightarrow e^{-tA}f\in L^2_+.$

\item $f\in \mathcal{D}(A^{\frac12})\Rightarrow |f|\in\mathcal{D}(A^{\frac12})\quad\hbox{and}\quad\|A^{\frac12}|f|\|_{L^2}\le\|A^{\frac12}f\|_{L^2}$
\end{enumerate} 
\end{Theo}
Now we will establish  the following result.
\begin{prop}\label{props23}
Let $d\in\{1,2,3\},\, p\in[1,\infty[,\alpha\geq0, \beta\in]0,1]$ and $ \lambda\geq e^{\frac{3+2\alpha}{\beta}}$. \mbox{Define $\mathcal{L}:=\frac{\DD^\beta}{\log^\alpha(\lambda+\DD)}$,} then 
\begin{enumerate}
\item
The family $(e^{-t\mathcal{L}})_{t\geq0}$ defines a $C_0$-semigroup of contractions in $L^p(\RR^d).
$ \item The family $(e^{-t\mathcal{L}})_{t\geq0}$ defines a sub-Markovian semigroup in $L^p(\RR^d).
$ 
\item The operator $(e^{-t\mathcal{L}})_{t\geq0}$ satisfies the Beurling-Deny theorem described in  Theorem $\ref{deurl}$.
\end{enumerate}
\end{prop}
\begin{proof}
{\bf($1$)} For $f\in L^p$ we define the action of the semigroup to this function by
$$
e^{-t\mathcal{L}}f(x)=K_t\star f(x),
$$
where $\widehat{K_t}(\xi)=e^{-t\frac{|\xi|}{\log^\alpha(\lambda+|\xi|)}}.$ From Corollary \ref{cor1} we have that the semigroup maps $L^p$ to itself  for every $p\in [1,\infty]$ and 
$$
\|K_t\star f\|_{L^p}\le\|f\|_{L^p}.
$$ The points $(1)$ and $(2)$ of the Definition \ref{def-sem}   are easy to check. It remains to prove the third point concerning the strong continuity of the semigroup. Since $\|K_t\|_{L^1}=1$ and $K_t\geq0,$ then  for  $\eta>0$ we have
\begin{eqnarray*}
K_t\star f(x)-f(x)&=&\int_{\RR^d}K_t(y)(f(x-y)-f(x))dy\\
&=&\int_{|y|\le\eta}K_t(y)(f(x-y)-f(x))dy\\
&+&\int_{|y|\geq\eta}K_t(y)(f(x-y)-f(x))dy\\
&=&\hbox{I}_1(x)+\hbox{I}_2(x).
\end{eqnarray*}
The first term is estimated as follows
\begin{eqnarray*}
\|\hbox{I}_1\|_{L^p}&\le&\int_{|y|\le\eta}K_t(y)\|f(\cdot-y)-f(\cdot)\|_{L^p}dy\\
&\le&\sup_{|y|\le\eta}\|f(\cdot-y)-f(\cdot)\|_{L^p}.
\end{eqnarray*}
For the second term we write
$$
\|\hbox{I}_2\|_{L^p}\le 2\|f\|_{L^p}\int_{|y|\geq \eta}K_t(y)dy.
$$
Combining these estimates we get
$$
\|K_t\star f-f\|_{L^p}\le\sup_{|y|\le\eta}\|f(\cdot-y)-f(\cdot)\|_{L^p}+2\|f\|_{L^p}\int_{|y|\geq \eta}K_t(y)dy.
$$
It is well-know that for every $p\in[1,\infty[$ we have
$$
\lim_{\eta\to 0^+}\sup_{|y|\le\eta}\|f(\cdot-y)-f(\cdot)\|_{L^p}=0.
$$
Thus for a given $\varepsilon>0$ we can find  $\eta> 0$ small enough such that
$$
\sup_{|y|\le\eta}\|f(\cdot-y)-f(\cdot)\|_{L^p}\le\varepsilon.
$$
Now to conclude the proof it suffices to prove that
$$
\lim_{t\to0^+}\int_{|y|\geq \eta}K_t(y)dy=0.
$$
This assertion is a consequence of the following result
$$
K_t\overset{t\to 0^+}{\rightharpoonup} \delta_0.
$$
To prove the last one we write for $\phi\in\mathscr{S},$
\begin{eqnarray*}
\langle K_t,\phi\rangle&=&\frac{1}{(2\pi)^d}\langle \widehat{K_t},\widehat\phi\rangle\\
&=&\frac{1}{(2\pi)^d}\int_{\RR^d}e^{-t|\xi|^\alpha}\widehat\phi(\xi)d\xi.
\end{eqnarray*}
We can use now Lebesgue theorem and the inversion Fourier transform leading to
$$
\lim_{t\to 0^+}\langle K_t,\phi\rangle=\frac{1}{(2\pi)^d}\int_{\RR^d}\widehat\phi(\xi)d\xi=\phi(0).
$$
Finally we get that $(K_t\star)_{t\geq0}$ defines a $C_0$-semigroup of contractions for every $p\in[1,\infty[$.
\vspace{1cm}

{\bf($2$)}
From the Definition \ref{defis2} and the first part of Proposition \ref{props23}  it remains to show that 
\begin{enumerate}
\item[$(1^\prime)$] For $f\in L^p$ with $  f(x)\geq 0 ,$ a.e. we have $e^{-t\mathcal{L}}f(x)\geq 0$
\item[$(2^\prime)$] For $f\in L^p$ with $  |f(x)|\leq 1$, a.e. we have $|e^{-t\mathcal{L}}f(x)|\le 1$
\end{enumerate}
The proof  is a direct consequence of the explicit formula
$$
e^{-t\mathcal{L}}f(x)=K_t\star f(x), 
$$
where according to Lemma \ref{kernel} we have $K_t\geq0$ and $\|K_t\|_{L^1}=1$.

{\bf{($3$)}} It is not hard to see that the operator $\frac{\DD^\beta}{\log^\alpha(\lambda+\DD)}$ is a nonnegative self-adjoint operator of $L^2$. This operator satisfies the first condition of Theorem \ref{deurl} since the kernel $K_t$ is positive.
\end{proof}

The following result gives in particular Theorem \ref{max-princ}.
\begin{prop}\label{maxi}Let $
A$ be a generator of  a $C_0$-semigroup of contractions, then
  \begin{enumerate}\item
Let $p\in[1,\infty[$ and $u\in \mathcal{D}({A})$.   then 
$$
\int_{\RR^2}A u\,|u|^{p-1}\textnormal{sign}\, u\,dx\leq 0.
$$
\item
Let $\theta$ be a smooth solution of the equation 
$$
\partial_t \theta+v\cdot\nabla\theta-A\theta=f
$$where $v$ is a smooth vector-field with zero divergence and $f$ a smooth function. Then for every $p\in[1,\infty]$ 
$$
\|\theta(t)\|_{L^p}\le\|\theta_0\|_{L^p}+\int_0^t\|f(\tau)\|_{L^p}d\tau.
$$

\end{enumerate}
\end{prop}
\begin{proof}
{\bf{(1)}} 
We introduce the operation $[h,g]$ between two functions by 
$$
[h,g]=\|g\|_{L^p}^{2-p}\int_{\RR^2}h(x)|g(x)|^{p-1}\textnormal{sign}\, g(x)dx.
$$
Now, we define the function  $\psi:[0,\infty[\to \RR$ by
$$
\psi(t)=\Big[e^{t A}u,u\Big].
$$
We have $\psi(0)=\|u\|_{L^p}^2$ and from H\"{o}lder inequality combined with the fact that the operator $e^{tA}$ is a contraction on $L^p$ we get
\begin{eqnarray*}
\psi(t)&\le&\|e^{t A}u\|_{L^p}\|u\|_{L^p}\\
&\le&\|u\|_{L^p}^2.
\end{eqnarray*}
Thus we find $\psi(t)\le\psi(0), \forall t\geq 0$. Therefore we get $\lim_{t\to0^+}\frac{\psi(t)-\psi(0)}{t}\leq 0.$ This gives
$$
\int_{\RR^2}Au(x)|u(x)|^{p-1}\textnormal{sign}\, u(x)dx\leq0.
$$

{\bf 2)}
Let $p\in[1,\infty[$ then multiplying the equation \eqref{transport-d} by $|\theta|^{p-1}\textnormal{sign} \theta$ and integrating by parts using $\textnormal{div} v=0$ we get
$$
\frac1p\frac{d}{dt}\|\theta(t)\|_{L^p}^p+\int_{\RR^2}|A\theta(x)\theta(x)|^{p-1}\textnormal{sign}\, \theta(x)dx\le\|f(t)\|_{L^p}\|\theta(t)\|_{L^p}^{p-1}.
$$
Using Proposition  \ref{maxi} we find
$$
\frac1p\frac{d}{dt}\|\theta(t)\|_{L^p}^p\le \|f(t)\|_{L^p}\|\theta(t)\|_{L^p}^{p-1}.
$$
By simplifying
$$
\frac{d}{dt}\|\theta(t)\|_{L^p}\le \|f(t)\|_{L^p}.
$$
Integrating in time we get for $p\in[1,\infty[$
$$
\|\theta(t)\|_{L^p}\le\|\theta_0\|_{L^p}.
$$
Since the estimates are uniform on the parameter $p$ then we can get the limit case $p=+\infty.$


\end{proof}

\subsection{Logarithmic estimate}
Let us now move to the last part of this section which deals with some logarithmic estimates generalizing  the results of \cite{vis,hk}. First we recall 
the following  result of propagation of Besov regularities.
\begin{prop}\label{prop-Bes}
Let  $\kappa\geq0$ and $A$ be a $C_0$ semigroup of contractions on $L^m(\RR^d)$ for every $m\in[1,\infty[.$ We assume that for every $q\in\NN\cup\{-1\}$, the operator $\Delta_q$ and $A$ commute on  a dense subset of $L^p$. Let $(p,r)\in[1,\infty]^2, s\in]-1,1[$ and $\theta$   be  a smooth solution of
$$
\partial_t\theta+v\cdot\nabla\theta-\kappa A\theta=f.
$$
  Then we have
$$
\|\theta\|_{\widetilde L^\infty_tB_{p,r}^s}\lesssim e^{CV(t)}\Big(\|\theta_0\|_{B_{p,r}^s}+\int_0^te^{-CV(\tau)}\|f(\tau)\|_{B_{p,r}^s}d\tau \Big),
$$
where $V(t)=\|\nabla v\|_{L^1_t L^\infty}$ and $C$  a constant depending only on $s$ and $d$.
\end{prop}
\begin{proof}
We set $\theta_q:=\Delta_q\theta$ then 
by localizing in frequency  the equation of $\theta$ we  get
$$
\partial_t\theta_q+v\cdot\nabla\theta_q-\kappa A\theta_q=-[\Delta_q, v\cdot\nabla]\theta+f_q.
$$
Using Proposition \ref{maxi} we get
$$
\|\theta_q(t)\|_{L^p}\le\|\theta_q(0)\|_{L^p}+\int_0^t\|[\Delta_q, v\cdot\nabla]\theta(\tau)\|_{L^p}d\tau+\int_0^t\|f_q(\tau)\|_{L^p}d\tau.
$$
On the other hand we have the classical  commuator estimate, see \cite{che1}
$$
\|[\Delta_q, v\cdot\nabla]\theta\|_{L^p}\le C 2^{-qs}c_q\|\nabla v\|_{L^\infty} \|\theta\|_{B_{p,r}^s},\quad \|(c_q)\|_{\ell^r}=1.
$$
Thus
$$
\|\theta(t)\|_{B_{p,r}^s}\le\|\theta_0\|_{B_{p,r}^s}+C\int_0^t\|\nabla v\|_{L^\infty} \|\theta\|_{B_{p,r}^s}+\int_0^t\|f(\tau)\|_{B_{p,r}^s}d\tau.
$$
It suffices now to use Gronwall inequality
\end{proof}

Now we will  show  that for the index regularity $s=0$ we can obtain a better estimate with a linear growth on the norm of the velocity. 

\begin{prop}
\label{thmlog}
Let $v$ be a smooth divergence free vector-field on $\RR^d$. Let  $\kappa\geq0$ and $A$ be a generator of $C_0$-semigroup  of contractions on $L^p(\RR^d)$ for every $p\in[1,\infty[.$ We assume that for every $q\in\NN$, the opeartors $\Delta_q$ and $A$ commute on  a dense subset of $L^p$. Let   $\theta$ be  a smooth solution of
$$
\partial_t\theta+v\cdot\nabla\theta-\kappa A\theta=f.
$$
Then we have for every $p\in [1,\infty]$
$$
\|\theta\|_{\widetilde L^\infty_tB_{p,1}^0}\leq C\Big( \|\theta_0\|_{B_{p,1}^0}+\|f\|_{L^1_tB_{p,1}^0}\Big)\Big( 1+\int_0^t\|\nabla v(\tau)\|_{L^\infty}d\tau \Big),
$$
where the constant $C$ does not depend on $p$ and $\kappa.$
\end{prop}

\begin{proof}
We mention that the result is first proved in \cite{vis} for the case $\kappa=0$ by using the special structure of the transport equation. In \cite{H-K2} Keraani and the author    generalized Vishik's result for a transport-diffusion equation where the dissipation term has the \mbox{form $-\kappa\Delta\theta$.}  The method described in \cite{H-K2} can be easily adapted  here for  our model. 

 Let 
$q\in\NN\cup\{-1\}$ and  denote by  $\overline\theta_q$  the unique global solution of the initial value problem
\begin{equation}\label{R_{Q}}\left\lbrace
\begin{array}{l}
\partial_t \overline\theta_q+v\cdot\nabla \overline\theta_q-\kappa A\overline\theta_q=\Delta_{q}f,
\\
{\overline\theta_q}_{|t=0}=\Delta_{q}\theta^{0}.\\
\end{array}
\right.
\end{equation}
Using  Proposition \ref{prop-Bes}  with $s=\pm\frac12$   we get
$$
\|\overline\theta_q\|_{\widetilde L^\infty_tB_{p,\infty}^{\pm\frac12}}\lesssim
\big(\|\Delta_{q} \theta_0\|_{B_{p,\infty}^{\pm\frac12}}+
\|\Delta_{q}f\|_{L^1_{t}B_{p,\infty}^{\pm\frac12}}\big) e^{CV(t)
},
$$
where $
V(t)=\|\nabla v\|_{L^1_{t}L^\infty}.
$
Combined with  the  definition of Besov spaces this yields \mbox{ for $j,q\geq-1$}
\begin{equation}
\label{t1}
\|\Delta_{j}\overline\theta_q\|_{L^\infty_tL^p}\lesssim 2^{-\frac12|j-q|}
\big(\|   \Delta_{q} \theta_0\|_{L^p}+\|\Delta_{q}f\|_{L^1_{t}L^p}\big) e^{CV(t)}.
\end{equation}
By linearity and again  the definition of Besov spaces we have
\begin{eqnarray}\label{t2}
\|\theta\|_{\widetilde L^\infty_tB_{p,1}^0}\leq
\sum_{|j-q|\geq N}
\|\Delta_{j}\overline\theta_q\|_{L^\infty_tL^p}+\sum_{|j-q|< N}
\|\Delta_{j}\overline\theta_q\|_{L^\infty_tL^p},
\end{eqnarray}
where $N\in \Bbb N$ is to be chosen later.
To deal with the first sum we use (\ref{t1}) 
\begin{eqnarray*}
\nonumber\sum_{|j-q|\geq N}
\|\Delta_{j}\overline\theta_q\|_{L^\infty_tL^p}&\lesssim& 2^{-N/2}\sum_{q\geq-1}\big(\|\Delta_{q}\theta_0\|_{L^p}+\|\Delta_{q}f\|_{L^1_{t}L^p}\big)e^{CV(t)}
\\
&\lesssim&  2^{-N/2}
\big(\|\theta ^0\|_{B_{p,1}^0}+\|f\|_{L^1_{t}B_{p,1}^0}\big)e^{CV(t)}.
\end{eqnarray*}
We now turn to the second  sum in  the right-hand side of (\ref{t2}).  

It is clear that
\begin{equation*}\nonumber\sum_{|j-q|< N}
\|\Delta_{j}\overline\theta_q\|_{L^\infty_tL^p}\lesssim \sum_{|j-q|< N}
\|\overline\theta_q\|_{L^\infty_tL^p}.
\end{equation*}
Applying Proposition \ref{propmax} to  the system (\ref{R_{Q}}) yields 
$$
\|\overline\theta_q\|_{L^\infty_tL^p}\leq \|\Delta_q\theta_0\|_{L^p}+\|\Delta_{q}f\|_{L^1_{t}L^p}.
$$
It follows that 
\begin{equation*}
\sum_{|j-q|<N}\|\Delta_{j}\overline\theta_q\|_{L^\infty_tL^p}\lesssim N\big(\|\theta ^0\|_{B_{p,1}^0}+\|f\|_{L^1_{t}B_{p,1}^0}\big).
\end{equation*}
The outcome is  the following
$$
\|\theta\|_{\widetilde L^\infty_tB_{p,1}^0}\lesssim\big(\|\theta ^0\|_{B_{p,1}^0}+\|f\|_{L^1_{t}B_{p,1}^0}\big)\Big(2^{- N/2}e^{CV(t)}+N\Big).
$$
Choosing
$$
N=\Big[\frac{2C
V(t)}{ \log 2}\Big]+1,
$$
we get the desired result.

\end{proof}
Combining Propositions \ref{thmlog} and  \ref{props23} we get,
\begin{cor}
\label{thmlog}
Let $v$ be a smooth divergence free vector-field on $\RR^d$, with $d\in\{2,3\}.$ Let  $\kappa,\alpha\geq 0,\beta\in]0,1],\,\lambda\geq e^{\frac{3+2\alpha}{\beta}}$, $ p\in[1,\infty]$ and  $\theta$ be  a smooth solution of
$$
\partial_t\theta+v\cdot\nabla\theta+\kappa\DD^\beta\log^{-\alpha}(\lambda+\DD)\theta=f.
$$
Then we have
$$
\|\theta\|_{\widetilde L^\infty_tB_{p,1}^0}\leq C\Big( \|\theta_0\|_{B_{p,1}^0}+\|f\|_{L^1_tB_{p,1}^0}\Big)\Big( 1+\int_0^t\|\nabla v(\tau)\|_{L^\infty}d\tau \Big),
$$
where the constant $C$ depends only on $\lambda$ and $\alpha$.
\end{cor}
\section{Proof of Theorem \ref{coer10}}
 
\subsection{Bernstein inequality}
This section is devoted to the generalization of  the classical Bernstein inequality described in \mbox{Lemma \ref{lb}} for more general operators.
\begin{prop}\label{Bernst} Let $\alpha\in\RR,\beta>0$ and $ \lambda\geq2.$ Then there exists a constant $C$ such that for every  $f\in\mathcal{S}(\RR^d)$ and  for every $q\geq-1$ and $ p\in[1,\infty]$ we have
$$
\Big\|\Delta_q \big(\frac{\DD^\beta}{\log^\alpha(\lambda+\DD)} f\big)\Big\|_{L^p}\le C2^{q\beta}{(|q|+1)^{-\alpha}} \|\Delta_qf\|_{L^p}.
$$
Moreover
$$
\Big\| S_q \big(\frac{\DD^\beta}{\log^\alpha(\lambda+\DD)} f\big)\Big\|_{L^p}\le C2^{q\beta}{(|q|+1)^{-\alpha}} \|S_qf\|_{L^p}.
$$
\end{prop}
\begin{Rema}\label{rmq6}
The first result of Proposition \ref{Bernst} remains true for more general situation where  $q\in\NN$ and the operator $\DD^\beta$ is replaced by $a(\textnormal{D})$ with $a(\xi)$ a homogeneous distribution of order $\beta\in\RR$ that is $a\in C^\infty(\RR\backslash\{0\})$  and for every $\gamma\in \NN^d$
$$
|\partial_\xi^\gamma a(\xi)|\le C|\xi|^{\beta-|\gamma|}.
$$
\end{Rema}
\begin{proof}

\underline{\it Case $q\in\NN$.}
It is easy to see that
$$
\Delta_q \big(\frac{\DD^\beta}{\log^\alpha(\lambda+\DD)} f\big)=K_q\star\Delta_qf,
$$
with 
$$
\widehat{K_q}(\xi)=\frac{\tilde{\phi}(2^{-q}\xi)|\xi|^{\beta}}{{\log^\alpha(\lambda+|\xi|)}}
$$
and $\tilde{\phi}$ is a smooth function   supported in the ring $\{\frac14\le|x|\le 3\}$ and taking the value 1 on the support of the function $\phi$ introduced in  section \ref{preliminaries}  . By Fourier inversion formula and change of variables we get
\begin{eqnarray*}
K_q(x)&=&c_d\int_{\RR^d}e^{ix\cdot\xi}\frac{\tilde{\phi}(2^{-q}\xi)|\xi|^{\beta}}{{\log^\alpha(\lambda+|\xi|)}}d\xi\\
&=&c_d2^{q\beta}2^{qd}\int_{\RR^d}e^{i2^qx\cdot\xi}\frac{\tilde{\phi}(\xi)|\xi|^{\beta}}{{\log^\alpha(\lambda+2^q|\xi|)}}d\xi\\
&:=&c_d2^{q\beta}2^{qd}\tilde{K_q}(2^q x),
\end{eqnarray*}
with
$$
\tilde{K_q}(x)=\int_{\RR^d}e^{ix\cdot\xi}\frac{\tilde{\phi}(\xi)|\xi|^{\beta}}{{\log^\alpha(\lambda+2^q|\xi|)}}d\xi.
$$
Obviously we have
$$
\|{K_q}\|_{L^1}=c_d2^{q\beta}\|\tilde{K_q}\|_{L^1}.
$$
Hence to prove Proposition \ref{Bernst} it suffices to establish
\begin{equation}\label{be1}
\|\tilde{K_q}\|_{L^1}\le C{(q+1)^{-\alpha}}.
\end{equation}
From the definition of $\tilde{K_q}$ we see that
$$
\tilde{K_q}(x)=\int_{\RR^d}e^{ix\cdot\xi}\frac{\tilde{\psi}(\xi)}{{\log^\alpha(\lambda+2^q|\xi|)}}d\xi
$$
where $\tilde\psi$ belongs to Schwartz class and supported in $\{\frac14\le|x|\le 3\}$. By integration by parts we get for $j\in\{1,2,...,d\}$
$$
x_j^{d+1}\tilde{K_q}(x)=(-i)^{d+1}\int_{\frac14\le|\xi|\le 3}e^{ix\cdot\xi}\,\partial_{\xi_j}^{d+1}\Big(\frac{\tilde{\psi}(\xi)}{{\log^\alpha(\lambda+2^q|\xi|)}}\Big)d\xi.
$$
Now we claim that
$$
\Big|\partial_{\xi_j}^{d+1}\Big(\frac{\tilde{\psi}(\xi)}{{\log^\alpha(\lambda+2^q|\xi|)}}\Big)  \Big|\le C_{\lambda,\alpha,d}\frac{g(\xi)}{{\log^\alpha(\lambda+2^{q})}},
$$
where $g\in \mathcal{S}(\RR^d).$ This is an easy consequence of Leibniz formula and the following fact
\begin{eqnarray*}
\Big|\partial_{\xi_j}^{n}\Big(\frac{1}{{\log^\alpha(\lambda+2^q|\xi|)}}\Big)\Big|&\le& \sum_{l,k=1}^nc_{l,k}\Big(\frac{2^q}{\lambda+2^q|\xi|}\big)^l\frac{1}{{\log^{\alpha+k}(\lambda+2^q|\xi|)}}\\
&\le&\frac{C_{\lambda,\alpha,n}}{{\log^{\alpha}(\lambda+2^{q})}},\quad\hbox{for}\quad \frac14\le|\xi|\le2.
\end{eqnarray*}
Thus we get for $j\in\{1,..,d\}$
$$
|x_j|^{d+1}|\tilde{K_q}(x)|\le {C}{{\log^{-\alpha}(\lambda+2^{q})}},\,\forall x\in\RR^d.
$$
It follows that
$$
|x|^{d+1}|\tilde{K_q}(x)|\le {C}{{\log^{-\alpha}(\lambda+2^{q})}},\,\forall x\in\RR^d.
$$
It is easy to see that  $\tilde{K_q}$ is continuous  and 
$$
|\tilde{K_q}(x)|\le{C}{{\log^{-\alpha}(\lambda+2^{q})}}
$$ 
Consequently,
$$
|\tilde{K_q}(x)|\le {C}{{\log^{-\alpha}(\lambda+2^{q})}}(1+|x|)^{-d-1},\,\forall x\in\RR^d.
$$
This yields
\begin{eqnarray*}
\|\tilde{K_q}\|_{L^1}&\le &{C}{{\log^{-\alpha}(\lambda+2^{q})}}\\
&\le& C(q+1)^{-\alpha}.
\end{eqnarray*}
This concludes the proof of the first case $q\in\NN$.

\underline{Case $q=-1$.} We can write in this case the kernel $K_{-1}$ as
\begin{eqnarray*}
K_{-1}(x)&=&\int_{\RR^d}e^{ix\cdot\xi}\frac{\tilde{\chi}(\xi)|\xi|^{\beta}}{{\log^\alpha(\lambda+|\xi|)}}d\xi\\
&=&\int_{\RR^d}e^{ix\cdot\xi}{{\chi_1}(\xi)}d\xi,
\end{eqnarray*}
where $\tilde\chi$ is  a smooth compactly supported function taking the value $1$ on the support of the function  $\chi$ introduced in section \ref{preliminaries}. The function $\chi_1$ is given by  ${\chi_1}(\xi)=\frac{\tilde{\chi}(\xi)|\xi|^{\beta}}{{\log^\alpha(\lambda+|\xi|)}}.$ We can see by easy computations  that $\tilde{\chi}$ is smooth outside zero and satisfies for every $\gamma\in\NN^d,$
$$
|\partial_\xi^\gamma\tilde{\chi}(\xi)|\le C_\gamma|\xi|^{\beta-|\gamma|},\quad\forall \xi \neq0.
$$
Using Mikhlin-H\"ormander theorem we get
$$
|K_{-1}(x)|\le C|x|^{-d-\beta}.
$$
Since $K_{-1}$ is continuous at zero  then we have $|K_{-1}(x)|\le C (1+|x|)^{-d-\beta}.$ This proves that $K_{-1}\in L^1$. 

To prove the second estimate we use the first result combined with the following  \mbox{identity $S_{q+2}S_q=S_q.$}
\begin{eqnarray*}
\Big\|S_q\big(\frac{\DD^\beta}{\log^\alpha(\lambda+\DD)} f\big)\Big\|_{L^p}&\le&\sum_{j=-1}^{q+1}\Big\|\Delta_j \big(\frac{\DD^\beta}{\log^\alpha(\lambda+\DD)} S_q f\big)\Big\|_{L^p}\\&\le& C\|S_q f\|_{L^p}\sum_{j=-1}^{q+1}2^{j\beta}(|j|+1)^{-\alpha}.
\end{eqnarray*}
Since $\beta>0$ then the last series diverges and 
$$
\sum_{j=-1}^{q+1}2^{j\beta}(|j|+1)^{-\alpha}\le C2^{q\beta} (|q|+1)^{-\alpha}.
$$
This can be deduced from the asymptotic behavior
$$
\int_{1}^{x} e^{\beta t}t^{-\alpha}dt\approx \frac1\beta e^{\beta x}x^{-\alpha},\quad\hbox{as}\quad x\to+\infty.
$$
\end{proof}
As a consequence of Proposition \ref{Bernst} we get the following result which describes  discuss the action of the logarithmic Riesz transform $\mathcal{R}_\alpha=\frac{\partial_1\log^{\alpha}(\lambda+\DD)}{\DD}$  on  Besov spaces.
\begin{cor}\label{embed}
Let $\alpha\in\RR,\lambda>1$ and $p\in[1,\infty]$.  Then the map
$$(\textnormal{Id}-\Delta_{-1})\mathcal{R}_\alpha :B_{p,r}^{s,\alpha}\to B_{p,r}^s
$$ is continuous.
\end{cor}

\subsection{Generalized Bernstein inequality}
The main goal of this section  is   to prove Theorem \ref{coer10}. Some preliminaries lemmas will be needed. The first one 
is a  Stroock-Varopoulos inequality for  sub-Markovian operators. For the proof see \cite{Lisk, Lisk1}.
\begin{Theo} Let $p>1$ and  $A$ be a sub-Markovian generator, then we have
$$
4\frac{p-1}{p^2}\|A^{\frac12}(|f|^{\frac{p}{2}}\textnormal{sign }f)\|_{L^2}^2\le\int_{\RR^d} (Af)\,|f|^{p-1}\textnormal{sign}\,fdx\leq C_p\|A^{\frac12}(|f|^{\frac{p}{2}}\textnormal{sign }f)\|_{L^2}^2.
$$
Moreover the generator $A$ satisfies the first Deurling-Deny condition
$$
4\frac{p-1}{p^2}\|A^{\frac12}(|f|^{\frac{p}{2}})\|_{L^2}^2\le\int_{\RR^d} (Af)\,|f|^{p-1}\textnormal{sign }fdx.
$$
\end{Theo}
Combining this result with Proposition \ref{props23} we get,
\begin{cor}\label{coer}
Let $, p>1, \beta\in]0,1],\alpha\geq0$ and $\lambda\geq e^{\frac{3+2\alpha}{\beta}}$. Then we have
$$
4\frac{p-1}{p^2}\Big\|\frac{\DD^{\frac\beta2}}{\log^{\frac\alpha2}(\lambda+\DD)}(|f|^{\frac{p}{2}})\Big\|_{L^2}^2\le\int_{\RR^d} \Big(\frac{\DD^{\beta}}{\log^{\alpha}(\lambda+\DD)}f\Big)\,|f|^{p-1}\textnormal{sign }fdx.
$$
\end{cor}

We will make use of  the  following composition results,
\begin{lem}\label{compos-bes}
\begin{enumerate}
\item Let $\mu\geq1$ and $s\in[0,\mu[\cap [0,2[.$ Then
$$
\||f|^\mu\|_{B_{2,2}^s}\le C\|f\|_{B_{2\mu,2}^s}\|f\|_{B_{2\mu,2}^0}^{\mu-1}
$$
\item $\mu\in]0,1], p,q\in[1,\infty]$ and $0<s<1+\frac1p.$ Then
$$
\||f|^\mu\|_{B_{\frac p\mu,\frac q\mu}^{s\mu}}\le C\|f\|_{B_{p,q}^s}^\mu.
$$
\end{enumerate}
\end{lem}
We point out that the first estimate is a particular case of a general result due to  Miao {\it et al.}, see \cite{cmz}. The second one is established by Sickel in \cite{Si}, see also Theorem 1.4 of \cite{kateb}.

Next we will recall the following result  proved in \cite{cmz,dan,plan},
\begin{prop}\label{coer01}
Let $d\geq1,\,\beta\in]0,2]$ and $p\geq2.$ Then we have for $q\in\NN$ and $f\in \mathcal{S}(\RR^d),$
$$
2^{q\beta}\|\Delta_q f\|_{L^p}^p\le C\int_{\RR^d} \big({\DD^{\beta}}\Delta_qf\big)\,|\Delta_q f|^{p-1}\textnormal{sign }\Delta_qfdx.
$$
where $C$ depends on $p$ and $\beta$.
Moreover, for $\beta=2$ we can extend the above inequality to $p\in]1,\infty[$.
\end{prop}
Now we will restate and prove Theorem \ref{coer10}.
\begin{prop}\label{coer1}
Let $d\in\{1,2,3\}, \beta\in]0,1],\alpha\geq0, \lambda\geq e^{\frac{3+2\alpha}{\beta}}$ and $p>1.$ Then we have for $q\in\NN$ and $f\in \mathcal{S}(\RR^d),$
$$
2^{q\beta}(q+1)^{-\alpha}\|\Delta_q f\|_{L^p}^p\le C\int_{\RR^d} \Big(\frac{\DD^{\beta}}{\log^{\alpha}(\lambda+\DD)}\Delta_qf\Big)\,|\Delta_q f|^{p-1}\textnormal{sign }\Delta_qfdx.
$$
where $C$ depends on $p,\alpha$ and $\lambda$.
\end{prop}
\begin{proof}
Using  Corollary \ref{coer} it suffices to prove
$$
C^{-1}2^{q\beta}(q+1)^{-\alpha}\|\Delta_q f\|_{L^p}^p\le \Big\|\frac{\DD^{\frac\beta2}}{\log^{\frac\alpha2}(\lambda+\DD)}(|\Delta_qf|^{\frac{p}{2}})\Big\|_{L^2}^2.
$$
We will use an idea of Miao {\it et al.} \cite{cmz}.
Let $N\in \NN$  then we have
$$
\|\DD(|f_q|^{\frac{p}{2}})\|_{L^2}\le\|S_N\DD(|f_q|^{\frac p2})\|_{L^2}+\|\|(\textnormal{Id}-S_N)\DD(|f_q|^{\frac p2})\|_{L^2}.
$$
It is clear that  for $s\geq0$
\begin{equation}\label{sick2}
\|\|(\textnormal{Id}-S_N)\DD(|f_q|^{\frac p2})\|_{L^2}\le C 2^{-N s}\||f_q|^{\frac p2}\|_{ B_{2,2}^{1+s}}.
\end{equation}
We have now to deal with fraction powers in Besov spaces. We will treat differently the cases $p> 2$ and $p\le2$.

{\underline{\it{Case $p>2.$}}} Combining   Lemma \ref{compos-bes}-(1)  with Bernstein inequality we get under the assumption  $0<s<\min(\frac p2-1,2),$
\begin{eqnarray*}
\||f_q|^{\frac p2}\|_{ B_{2,2}^{1+s}}&\le &C\|f_q\|_{B_{p,2}^{0}}^{\frac p2-1}\|f_q\|_{B_{p,2}^{1+s}}\\
&\le& C2^{q(1+s)}\|f_q\|_{L^p}^{\frac p2}.
\end{eqnarray*}
\underline{\it{Case $1<p\le2.$}} Using Lemma \ref{compos-bes}-(2) and Bernstein inequality, we get for $0<s<\frac{p-1}{2},$
\begin{eqnarray*}
\||f_q|^{\frac p2}\|_{ B_{2,2}^{1+s}}&\le &C\|f_q\|_{B_{p,p}^{\frac{2+2s}{p}}}^{\frac p2}\\
&\le& C2^{q(1+s)}\|f_q\|_{L^p}^{\frac p2}.
\end{eqnarray*}

It follows from (\ref{sick2}) and the previous inequalities that there exists $s_p>0$ such that \mbox{for $0<s<s_p$}
$$
\|\|(\textnormal{Id}-S_N)\DD(|f_q|^{\frac p2})\|_{L^2}\le C2^{-N s}2^{q(1+s)}\|f_q\|_{L^p}^{\frac p2}.
$$
On the other hand Proposition \ref{Bernst} gives
\begin{eqnarray*}
\|S_N\DD(|f_q|^{\frac p2})\|_{L^2}&\le &\Big\|S_N{\DD^{1-\frac\beta2}}{\log^{\frac\alpha2}(\lambda+\DD)}\big(\frac{\DD^{\frac\beta2}}{\log^{\frac\alpha2}(\lambda+\DD)}(|f_q|^{\frac p2})\big)\Big\|_{L^2}\\
&\le&C 2^{N(1-\frac\beta2)}N^{\frac\alpha2}\Big\|\frac{\DD^{\frac\beta2}}{\log^{\frac\alpha2}(\lambda+\DD)}(|f_q|^{\frac p2})\big)\Big\|_{L^2}.
\end{eqnarray*}
Therefore we get
$$
\|\DD(|f_q|^{\frac{p}{2}})\|_{L^2}\le C2^{-N s}2^{q(1+s)}\|f_q\|_{L^p}^{\frac p2}+C 2^{N(1-\frac\beta2)}N^{\frac\alpha2}\Big\|\frac{\DD^{\frac\beta2}}{\log^{\frac\alpha2}(\lambda+\DD)}(|f_q|^{\frac p2})\big)\Big\|_{L^2}.
$$
According to Proposition \ref{coer01} we have for $p\in]1,\infty[$
$$C_p2^q\|f_q\|_{L^p}^{\frac p2}\le
\|\DD(|f_q|^{\frac{p}{2}})\|_{L^2}.
$$
Combining both last estimates we get
$$
2^q\|f_q\|_{L^p}^{\frac p2}\le C2^{s(q-N)}\,2^{q}\|f_q\|_{L^p}^{\frac p2}+C 2^{N(1-\frac\beta2)}N^{\frac\alpha2}\Big\|\frac{\DD^{\frac\beta2}}{\log^{\frac\alpha2}(\lambda+\DD)}(|f_q|^{\frac p2})\big)\Big\|_{L^2}
$$
We take $N=q+N_0$ such that $C2^{-N_0 s}\le\frac12$ . Then we get
$$
\|f_q\|_{L^p}^{\frac p2}\le C 2^{-q\frac\beta2}(q+1)^{\frac\alpha2}\Big\|\frac{\DD^{\frac\beta2}}{\log^{\frac\alpha2}(\lambda+\DD)}\big(|f_q|^{\frac p2}\big)\Big\|_{L^2}.
$$
This gives the desired result.
\end{proof}
\section{Commutator estimates}
We will establish in this section some commutator estimates. The following result was proved in \cite{HKR}.
\begin{lem}
 \label{commu}Given $(p,m)\in[1,\infty]^2$ such that $p\geq m'$ with $m'$ the conjugate exponent of $m$. Let $f,g$ and $h$ be three functions such that $\nabla f\in L^p, g\in L^m$ and $xh\in L^{m'}$. Then,
 $$
 \|h\star(fg)-f(h\star g)\|_{L^p}\leq \|xh\|_{L^{m'}}\|\nabla f\|_{L^p}\|g\|_{L^{m}}.
 $$

 \end{lem}
 Now we will prove the following lemma.
\begin{lem}\label{conv}
Let $(a_n)_{n\in\ZZ}$ be a sequence of strictly nonnegative real numbers such that 
$$
M:=\max\Big(\sup_{n\in\ZZ} a_n^{-1}\sum_{j\le n}a_j,\sup_{n\in\ZZ} a_n\sum_{j\geq n}a_j^{-1}\Big)<\infty.
$$
Then for every $p\in[1,\infty]$  the linear operator $T:\ell^p\to \ell^p$ defined by
$$
T((b_n)_{n\in\ZZ})=\Big(\sum_{j\le n}a_j a_n^{-1} b_j\Big)_{n\in\ZZ}
$$
is continuous and $\|T\|_{\mathcal{L}(\ell^p)}\le M.$
\end{lem}
 \begin{proof}
 By interpolation it suffices to prove the cases $p=1$ and $p=+\infty$. Let's start \mbox{with $p=1$} and denote ${\bf{b}}=(b_n)_{n\in\ZZ}.$ Then from Fubini lemma and the hypothesis
 \begin{eqnarray*}
 \|T{\bf b}\|_{\ell^1}&\le& \sum_{n\in\ZZ}\sum_{j\le n}a_j a_n^{-1} |b_j|\\
 &\le& \sum_{j\in\ZZ}|b_j|a_j\sum_{n\geq j}a_n^{-1} \\
 &\le& M\|{\bf b}\|_{\ell^1}.
 \end{eqnarray*}
 For the case $p=+\infty$ we write
 \begin{eqnarray*}
 \|T{\bf b}\|_{\ell^\infty}&\le& \sup_{n\in\ZZ}\sum_{j\le n}a_j a_n^{-1} |b_j|\\
 &\le& \|{\bf b}\|_{\ell^\infty}\sup_{n\in\ZZ}a_n^{-1}\sum_{j\leq n}a_j \\
 &\le& M\|{\bf b}\|_{\ell^\infty}.
 \end{eqnarray*}
 This completes the proof.
 \end{proof}


The goal now is to study the commutation between the following operators
$$\mathcal{R}_\alpha=\frac{\partial_1}{\DD}\log^\alpha(\lambda+\DD)\quad\hbox{and}\quad v \cdot \nabla. 
$$  Recall that $B_{\infty,2}^{s,s'}$ is the space given by the set of tempered distributions $u$ such that
 $$
 \|u\|_{B_{\infty,r}^{s,s'}}=\|\big(2^{qs}(|q|+1)^{s'}\|\Delta_qu\|_{L^\infty}\big)_q\|_{\ell^r}.
 $$
 The main result of this section reads as follows.
\begin{prop}\label{propcom}
Let $\alpha\in\RR,\lambda>1$, $v$ be   a smooth  divergence free  vector-field and $\theta $ be a smooth scalar function.  
\begin{enumerate}
\item For every $(p,r)\in [2,\infty[\times[1,\infty]$ there exists 
a constant $C=C(p,r)$ such that 
$$
\|[\mathcal{R_\alpha}, v\cdot\nabla]\theta\|_{B_{p,r}^0}\le C \|\nabla v\|_{L^p}\big(\|\theta\|_{B_{\infty,r}^{0,\alpha}}+\|\theta\|_{L^p}\big).
$$

\item 
For every  $(r,\rho)\in[1,\infty]\times]1,\infty[$  and $\epsilon >0$ there exists 
a constant $C=C(r,\rho,\epsilon)$ such that 

 $$
 \|[\mathcal{R_\alpha}, v\cdot\nabla]\theta\|_{B_{\infty,r}^0}\le C (\|\omega\|_{L^\infty}+\|\omega\|_{L^\rho})\big(\|\theta\|_{B_{\infty,r}^\epsilon}+\|\theta\|_{L^\rho}\big).
 $$
 \end{enumerate}
 \end{prop}
 \begin{proof}
 {\bf $(1)$}
 We  split the commutator into three parts according to  Bony's decomposition \cite{b},
\begin{eqnarray*}
\nonumber[\mathcal{R}_\alpha, v\cdot\nabla]\theta&=&\sum_{q\in\NN}[\mathcal{R}_\alpha, S_{q-1}v\cdot\nabla]\Delta_q\theta+\sum_{q\in\NN}[\mathcal{R}_\alpha, \Delta_qv\cdot\nabla]S_{q-1}\theta\\
\nonumber&+&\sum_{q\geq-1} [\mathcal{R}_\alpha, \Delta_qv\cdot\nabla]\widetilde{\Delta}_q\theta\\
\nonumber&=& \sum_{q\in\NN}\mbox{I}_q+\sum_{q\in\NN}\mbox{II}_q+\sum_{q\geq-1}\mbox{III}_q\\
&=&\mbox{I}+\mbox{II}+\mbox{III}.
\end{eqnarray*}
We start with the estimate of the first term $\mbox{I}$.    It is easy to see that there exists $\tilde\varphi\in\mathcal{S}$  whose spectum does not meet the origin such that 
 $$
\mbox{I}_q(x)=h_q\star( S_{q-1}v\cdot\nabla\Delta_q\theta)-S_{q-1}v\cdot(h_q\star \nabla\Delta_q\theta),
$$
 where $$
 \widehat{h_q}(\xi)=i\tilde\varphi(2^{-q}\xi)\frac{\xi_1}{|\xi|}\log^\alpha(\lambda+|\xi|).
 $$
  Applying Lemma \ref{commu} with $m=\infty$   we get 
  \begin{eqnarray}
  \nonumber
\|\mbox{I}_q\|_{L^p}&\lesssim &\|xh_q\|_{L^1}  \|\nabla S_{q-1} v\|_{L^p}\|\Delta_q\nabla\theta\|_{L^\infty}
\\
\label{x1}
&\lesssim& 2^q\|xh_q\|_{L^1}\|\Delta_q\theta\|_{L^\infty}\|\nabla v\|_{L^p}.
\end{eqnarray}
We can easily check that
$$
\|xh_q\|_{L^1}=  2^{-q}\|x\tilde{h}_q\|_{L^1}\quad\hbox{with}\quad   \widehat{\tilde{h}_q}(\xi)=i\tilde\varphi(\xi)\frac{\xi_1}{|\xi|}\log^\alpha(\lambda+2^q|\xi|).
$$
We can get by a similar way to the proof of  Proposition \ref{Bernst} 
$$
\|\tilde h_q\|_{L^1}\le C (1+|q|)^\alpha.
$$
Thus estimate (\ref{x1}) becomes
$$
\|\mbox{I}_q\|_{L^p}\le C (1+|q|)^\alpha\|\Delta_q\theta\|_{L^\infty}\|\nabla v\|_{L^p}.
$$
 Combined with the trivial fact 
 $$\Delta_j\sum_{q}\mbox{I}_q= \sum_{|j-q|\le 4}\mbox{I}_q
 $$ this yields
\begin{eqnarray*}
\|\mbox{I}\|_{B_{p,r}^0}
&\lesssim&
\Big(\sum_{q\geq-1}\|\mbox{I}_q\|_{L^p}^r\Big)^{\frac1r}
\\
&\lesssim&\|\nabla v\|_{L^p}
\|\theta\|_{B_{\infty,r}^{0,\alpha}}.
\end{eqnarray*}

Let us move to  the second term $\mbox{II}$. As before one writes
$$
\mbox{II}_q(x)=h_q\star( \Delta_q v\cdot\nabla  S_{q-1}\theta)-\Delta_q v\cdot(h_q\star \nabla S_{q-1}\theta),
$$ 
and then we obtain the estimate 
\begin{eqnarray*}
\nonumber\|\mbox{II}_q\|_{L^p}&\lesssim& 2^{-q}(1+|q|)^\alpha\|\Delta_q \nabla v\|_{L^p}\| S_{q-1}\nabla\theta\|_{L^\infty}\\
&\lesssim&\|\nabla v\|_{L^p}\sum_{j\le q-2}\frac{2^{j}(1+|j|)^{-\alpha}}{2^{q}(1+|q|)^{-\alpha}}((1+|j|)^\alpha\|\Delta_j\theta\|_{L^\infty}).
\end{eqnarray*}
Combined with Lemma \ref{conv} this  yields
$$
\|\mbox{II}\|_{B_{p,r}^0}\lesssim \|\nabla v\|_{L^p}\|\theta\|_{B_{\infty,r}^{0,\alpha}}.
$$

\

Let us now deal with the third term  $\mbox{III}$. Using  the fact that the divergence of $\Delta_q v$ vanishes, then  we can write  $\mbox{III}$  as
\begin{eqnarray*}
\mbox{III}&=&\sum_{q\geq 2} \mathcal{R}_\alpha \textnormal{div} (\Delta_qv\, \widetilde{\Delta}_q\theta)- \sum_{q\geq 2} \textnormal{div}(\Delta_qv\,\mathcal{R}_\alpha\widetilde{\Delta}_q\theta)+\sum_{q\leq 1}[\mathcal{R}_\alpha, \Delta_{q}v\cdot\nabla]\widetilde{\Delta}_{q}\theta\\&=&J_1+J_2+J_3.
\end{eqnarray*}
Using   Remark \ref{rmq6} we get 
\begin{eqnarray*}
\big\|\Delta_j\mathcal{R}_\alpha\textnormal{div} (\Delta_qv\, \widetilde{\Delta}_q\theta)\big\|_{L^p}\lesssim  2^j(1+|j|)^\alpha\|\Delta_q v\|_{L^p}\|\widetilde{\Delta}_q\theta\|_{L^\infty}.
\end{eqnarray*}
and since $q\geq2$
\begin{eqnarray*}
 \big\|\Delta_j \textnormal{div}(\Delta_qv\, \mathcal{R}_\alpha\widetilde{\Delta}_q\theta)\big\|_{L^p}&\lesssim & 2^j\|\Delta_q v\|_{L^p}\|\mathcal{R}_\alpha\widetilde{\Delta}_q\theta\|_{L^\infty}
 \\
 &\lesssim& 2^j(1+|q|)^{\alpha}\|\Delta_q v\|_{L^p}\|\widetilde{\Delta}_q\theta\|_{L^\infty}.
\end{eqnarray*}
Therefore we get
\begin{eqnarray*}
\|\Delta_j (J_1+J_2)\|_{L^p}&\lesssim&  \sum_{q\in\NN\atop q\geq j-4}2^j(1+|q|)^\alpha\|\Delta_q v\|_{L^p}\|\widetilde{\Delta}_q\theta\|_{L^\infty}\\
&\lesssim&\|\nabla v\|_{L^p} \sum_{q\in\NN\atop q\geq j-4}2^{j-q}(1+|q|)^\alpha\|{\Delta}_q\theta\|_{L^\infty},
\end{eqnarray*}
where we have again used Bernstein inequality to get the last line.
It suffices now to use Lemma \ref{conv}
$$
\|J_1+J_2\|_{B_{p,r}^0}\lesssim \|\nabla v\|_{L^p}\|\theta\|_{B_{\infty,r}^{0,\alpha}}.
$$
For the last term $J_3$ we can  write
$$
\sum_{-1\leq q\leq 1}[\mathcal{R}_\alpha, \Delta_{q}v\cdot\nabla]\widetilde{\Delta}_{q}\theta(x)=\sum_{q\leq 1}[\textnormal{div }\widetilde{\chi}(\textnormal{D})\mathcal{R}_\alpha, \Delta_{q}v]\widetilde{\Delta}_{q}\theta(x),$$
where  $\widetilde{\chi}$ belongs to $\mathcal{D}(\RR^d)$.
From the proof of Proposition  \ref{Bernst} we get that 
$\textnormal{div }\widetilde{\chi}(\textnormal{D})\mathcal{R}_\alpha$ is a convolution operator  with  a kernel   $\tilde{h}$ satisfying 
$$
|\tilde{h}(x)|\lesssim (1+|x|)^{-d-1}.
$$  
Thus
$$
J_3=  \sum_{q\leq 1} \tilde h\star( \Delta_qv\cdot\tilde\Delta_q\theta)-\Delta_qv\cdot(\tilde h\star \tilde\Delta_q\theta).
$$
First of all we point out that $\Delta_j J_3=0$ for $j\geq 6$, thus we just need  to estimate the low frequencies of $J_3$.   Noticing that $x\tilde h$ belongs to $L^{p'}$ for  $p'>1$ then using Lemma \ref{commu} with $m=p\geq 2$ we obtain
\begin{eqnarray*}
\|\Delta_j J_3\|_{L^p}&\lesssim&  \sum_{q\le 1}  \|x\tilde h\|_{L^{p'}} \|\Delta_q \nabla v\|_{L^p}\|\widetilde{\Delta}_q\theta\|_{L^p}\\
&\lesssim& \| \nabla v\|_{L^p}\sum_{-1\leq q\leq 1}\|{\Delta}_q\theta\|_{L^p}.
\end{eqnarray*}
This yields  finally
$$
\|J_3\|_{B_{p,r}^0}\lesssim \|\nabla v\|_{L^p}\|\theta\|_{L^p}.
$$
This completes the proof of  the first part of Theorem \ref{propcom}. 

{\bf$(2)$} The second part can be done in the same way so we will  give here just a shorten  proof. To estimate the terms ${\rm I}$ and ${\rm II}$ we use two facts: the first one is $\|\Delta_q\nabla u\|_{L^\infty}\approx \|\Delta_q \omega\|_{L^\infty}$ for all  $q\in \mathbb N.$ The second one is 
 \begin{eqnarray*}
\|\nabla S_{q-1} v\|_{L^\infty}& \lesssim&\|\nabla\Delta_{-1}v\|_{L^\infty}+\sum_{j=0}^{q-2}\|\Delta_j\nabla v\|_{L^\infty}\\
&\lesssim&\|\omega\|_{L^\rho}+ q\|\omega\|_{L^\infty}.
\end{eqnarray*}
Thus \eqref{x1} becomes
$$
\|\hbox{I}_q\|_{L^\infty}\le \|\omega\|_{L^\infty}(1+|q|)^{1+\alpha}\|\Delta_q\theta\|_{L^\infty}
$$
and by Proposition \ref{embed}
\begin{eqnarray*}
\|\hbox{I}\|_{B_{\infty,r}^0}&\le&\|\omega\|_{L^\infty}\|\theta\|_{B_{\infty,r}^{0,1+\alpha}}\\
&\le&\|\omega\|_{L^\infty}\|\theta\|_{B_{\infty,\infty}^\epsilon}.
\end{eqnarray*}
The second term $\hbox{II}$ is estimated as follows
\begin{eqnarray*}
\|\hbox{II}\|_{B_{\infty,r}^0}&\le&\|\omega\|_{L^\infty}\|\theta\|_{B_{\infty,r}^{0,\alpha}}\\
&\le&\|\omega\|_{L^\infty}\|\theta\|_{B_{\infty,\infty}^\epsilon}
\end{eqnarray*}

For the remainder term we do strictly the same analysis as before except for $J_3$: we apply Lemma \ref{commu} with $p=\infty$ and $m=\rho$ leading to 
\begin{eqnarray*}
\nonumber
\|\Delta_j J_3\|_{L^p}&\lesssim&  \sum_{q\le 1}  \|x\tilde h\|_{L^{\rho'}} \|\Delta_q \nabla v\|_{L^\infty}\|\widetilde{\Delta}_q\theta\|_{L^\rho}\\
&\lesssim& \| \nabla v\|_{L^\rho}\sum_{-1\leq q\leq 1}\|{\Delta}_q\theta\|_{L^\rho}\\
&\lesssim& \|\omega\|_{L^\rho}\|\theta\|_{L^\rho}.
\end{eqnarray*}
This ends  the proof of the theorem.

\end{proof}
\section{Smoothing effects}
In this section we will describe some smoothing effects for the model (\ref{transport-d}) and focus only on the case $\beta=1.$ Remark that we can obtain similar results for the case $\beta\in]0,1].$
\begin{equation} 
\left\{ \begin{array}{ll} 
\partial_{t} \theta+v\cdot\nabla \theta+\frac{{\DD}}{\log^\alpha(\lambda+\DD)}\vert \theta=f\\
\theta_{| t=0}=\theta^{0}.
\end{array} \right. \tag{${\textnormal{TD}}$}
\end{equation} 
%
We intend to prove the following smoothing effect. 
\begin{Theo}
 \label{thm99} Let  $\alpha\geq 0, \lambda\geq e^{3+2\alpha}, d\in\{2,3\}$ and $v$ be a smooth divergence-free vector field of $\RR^d$ with vorticity $\omega$. Then, for every   $p\in ]1,\infty[$  there exists a constant $C$ such that 
$$
\sup_{q\in\NN}2^q(1+q)^{-\alpha}\|\Delta_q \theta\|_{L^1_tL^p}\leq C\| \theta_0\|_{L^p}+C\|\theta_0\|_{L^\infty} \|\omega\|_{L^1_tL^p},
$$
for every  smooth  solution $\theta$ of ${\rm (TD)}$ with zero source term $f$.
\end{Theo}
\begin{Rema}
For the sake of simplicity we state the result of smoothing effect only for $\beta=1$ but the result remains true under the hypothesis of Proposition \ref{coer1}. 
\end{Rema}
\begin{proof}
We start with localizing in frequencies the equation:
for $q\geq-1$ we set $\theta_q:=\Delta_q\theta. $ Then 
$$
\partial_t\theta_q+v\cdot\nabla\theta_q+\frac{{\DD}}{\log^{\alpha}(\lambda+\DD)}\theta_q=-[\Delta_q, v\cdot\nabla]\theta.
$$
Recall that ${\theta}_q$ is real function since the functions involved in the dyadic partition of the unity are radial. Then multiplying the above equation  by $|\theta_q|^{p-2}{\theta}_q,$ integrating by parts  and using H\"{o}lder inequalities we get
$$
\frac1p\frac{d}{dt}\|\theta_q\|_{L^p}^p+\int_{\RR^2}\Big(\frac{\vert\textnormal{D}\vert}{\log^{\alpha}(\lambda+\DD)}\theta_q\Big) |\theta_q|^{p-2}{\theta}_qdx\leq \|\theta_q\|_{L^p}^{p-1}\|[\Delta_q, v\cdot\nabla]\theta\|_{L^p}.
$$
Using Proposition \ref{coer1} we get for $q\geq0$
$$
c  2^{q}(1+q)^{-\alpha}\|\theta_q\|_{L^p}^p\le\int_{\RR^2}\Big(\frac{\vert\textnormal{D}\vert}{\log^{\alpha}(\lambda+\DD)}\theta_q\Big) |\theta_q|^{p-2}{\theta}_qdx,
$$
where $c$ depends on $p$. Inserting this estimate in the previous one  we obtain 
$$
\frac1p\frac{d}{dt}\|\theta_q\|_{L^p}^p+c2^q(1+q)^{-\alpha}\|\theta_q\|_{L^p}^p\lesssim  \|\theta_q\|_{L^p}^{p-1} \|[\Delta_q, v\cdot\nabla]\theta\|_{L^p}.
$$
Thus we find 
\begin{equation}
\label{est}
\frac{d}{dt}\|\theta_q\|_{L^p}+c2^q(1+q)^{-\alpha} \|\theta_q\|_{L^p}\lesssim  \|[\Delta_q, v\cdot\nabla]\theta\|_{L^p}.
\end{equation}
To estimate the right hand-side we will use the following result,  see Proposition 3.3 of \cite{HKR}.
$$
\|[\Delta_q, v\cdot\nabla]\theta\|_{L^p}\lesssim \|\nabla v\|_{L^p}\|\theta\|_{B_{\infty,\infty}^{0}}.
$$

Combined with \eqref{est} this lemma yields
\begin{eqnarray*}
\frac{d}{dt}\big( e^{ct2^q(1+q)^{-\alpha}}\|\theta_q(t)\|_{L^p}\big)&\lesssim&   e^{ct2^q(1+q)^{-\alpha}} \|\nabla v(t)\|_{L^p}\|\theta(t)\|_{B_{\infty,\infty}^0}\\
&\lesssim& e^{ct2^q(1+q)^{-\alpha}}  \|\omega(t)\|_{L^p}\|\theta_0\|_{L^\infty}.
\end{eqnarray*}
To get the last line, we have used   the conservation of the $L^\infty$ norm of $\theta$ and the classical fact 
$$
 \|\nabla v \|_{L^p } \lesssim \| \omega \|_{L^p}\qquad\forall p\! \in ]1, +\infty[.
$$
Integrating the differential inequality  we get for $q\in\NN$
\begin{eqnarray*}
\|\theta_q(t)\|_{L^p}&\lesssim& \|\theta_q^0\|_{L^p} e^{-ct2^q(1+q)^{-\alpha}}+\|\theta_0\|_{L^\infty}\int_0^t e^{-c(t-\tau)2^q(1+q)^{-\alpha}} \|\omega(\tau)\|_{L^p}d\tau.
\end{eqnarray*}
Integrating in time yields finally
\begin{eqnarray*}
2^q(1+q)^{-\alpha}\|\theta_q\|_{L^1_tL^p}&\lesssim& \|\theta_q^0\|_{L^p} +\|\theta_0\|_{L^\infty}\int_0^t  \|\omega(\tau)\|_{L^p}d\tau\\
&\lesssim&\|\theta_0\|_{L^p} +\|\theta_0\|_{L^\infty}\int_0^t  \|\omega(\tau)\|_{L^p}d\tau,
\end{eqnarray*}
which is the desired result.
\end{proof}


\section{Proof of Theorem \ref{theo1} }

Throughout this section we use the notation $\Phi_k$ to denote any function
of the form 
$$
\Phi_k(t)=  C_{0}\underbrace{ \exp(...\exp  }_{k\,times}(C_0t)...),
$$
where $C_{0}$ depends on the involved norms of the initial data and its value may vary from line to line up to some absolute constants. 
We will make an intensive  use (without mentionning it) of  the following trivial facts
$$
\int_0^t\Phi_k(\tau)d\tau\leq \Phi_k(t)\qquad{\rm and}\qquad \exp({\int_0^t\Phi_k(\tau)d\tau})\leq \Phi_{k+1}(t).
$$

\label{sectionbouss}
The proof of Theorem \ref{theo1} will be done in several steps. The first one deals with  some {\it {\it a priori }   } estimates for the equations  \eqref{Bouss}. In the  second one we prove the uniqueness part. Finally,   we will  discuss    the construction of the solutions at the end of this section.
\subsection{{A priori    estimates}}

In Theorem \ref{theo1} we  deal with critical regularities  and   one needs to bound the Lipschitz norm of the velocity in order to get the global persistence of the initial regularities. For this purpose  we will proceed in several steps: one of the main steps is to  give an $L^\infty$-bound of the vorticity but due to some technical difficulties related to Riesz transforms  this will  not be done in a straight way. We prove before  an $L^p$ estimate for the vorticity with $2<p<\infty$. 

\subsubsection{ $L^p$-estimate of the vorticity}
We intend  now  to bound  the $L^p$-norm of the vorticity and to describe a smoothing effect for the temperature.
\begin{prop}\label{max-pro}
Let $\alpha\in[0,\frac12], \lambda\geq e^{3+2\alpha}$ and $p\in]2,\infty[$. Let $(v,\theta)$ be a solution of \eqref{Bouss} with   $\omega^0\in L^p,\, \theta_0\in L^p\cap L^\infty$ and $ \mathcal{R}_\alpha\theta_0\in L^p$.  Then for every $\epsilon>0$
$$
\|\omega(t)\|_{L^p}+\|\theta\|_{L^1_t B_{p,1}^{1-\epsilon}}\le \Phi_2(t).
$$
\end{prop}
\begin{proof}
 Applying the  transform $\mathcal{R}_\alpha$ to the temperature equation  we get
\begin{equation}\label{bm}
\partial_t\mathcal{R}_\alpha\theta+v\cdot\nabla\mathcal{R}_\alpha\theta+\frac{|\textnormal{D}|}{\log^\alpha(\lambda+\DD)}\mathcal{R}_\alpha\theta=-[\mathcal{R}_\alpha, v\cdot\nabla]\theta.
\end{equation}
Since $\frac{|\textnormal{D}|}{\log^\alpha(\lambda+\DD)}\mathcal{R}_\alpha=\partial_1,$ then the function $\Gamma:=\omega+\mathcal{R}_\alpha\theta$ satisfies 
\begin{equation}
\label{bm1}
\partial_t\Gamma+v\cdot\nabla\Gamma=-[\mathcal{R}_\alpha, v\cdot\nabla]\theta.
\end{equation}
According to  the first part of Proposition \ref{propcom} applied with $r=2$, 
$$
\big\|[\mathcal{R}_\alpha,v\cdot\nabla]\theta\big\|_{B_{p,2}^{0}}\lesssim\|\nabla v\|_{L^p}\big(\|\theta\|_{B_{\infty,2}^{0,\alpha}}+\|\theta\|_{L^p} \big).
$$
Using the classical   embedding $B_{p,2}^0\hookrightarrow L^p$ which is  true only for $p\in[2,\infty)$  
$$
\big\|[\mathcal{R}_\alpha,v\cdot\nabla]\theta\big\|_{L^p}\leq\|\nabla v\|_{L^p}\big(\|\theta\|_{B_{\infty,2}^{0,\alpha}}+\|\theta\|_{L^p}\big).
$$
Since ${\rm div }\,v=0$ then the $L^p$ estimate applied to   the transport equation \eqref{bm1} gives
$$
\|\Gamma(t)\|_{L^p}\leq\|\Gamma^0\|_{L^p}+\int_0^t\|[\mathcal{R}_\alpha, v\cdot\nabla]\theta(\tau)\|_{L^p}d\tau.
$$
Applying Theorem \ref{max-princ} to (\ref{bm}) yields
$$
\|\mathcal{R}_\alpha\theta(t)\|_{L^p}\le\|\mathcal{R}_\alpha\theta_0\|_{L^p}+\int_0^t\|[\mathcal{R}_\alpha, v\cdot\nabla]\theta(\tau)\|_{L^p}d\tau.
$$
We set $f(t):=\|\omega(t)\|_{L^p}+\|\mathcal{R}_\alpha\theta(t)\|_{L^p}$. Then from the previous  estimates  we get
\begin{eqnarray*}
f(t)&\lesssim& \|\Gamma_0\|_{L^p}+\|\mathcal{R}_\alpha\theta_0\|_{L^p}+\int_0^t\|\nabla v(\tau)\|_{L^p}\big(\|\theta(\tau)\|_{B_{\infty,2}^{0,\alpha}}+\|\theta\|_{L^p}\big) d\tau
\\
&\lesssim&f(0)+\int_0^t f(\tau)\big(\|\theta(\tau)\|_{B_{\infty,2}^{0,\alpha}}+\|\theta_0\|_{L^p}\big) d\tau.
\end{eqnarray*}
We have used here two estimates: the Calder\'on-Zygmund estimate: for $p\in(1,\infty)$
$$
\|\nabla v\|_{L^p}\le C\|\omega\|_{L^p}.
$$
The second one is $\|\theta(t)\|_{L^p}\le\|\theta_0\|_{L^p}$ described in Theorem \ref{max-princ}.

According to Gronwall lemma we get
\begin{equation}
\label{gr1}
f(t)\lesssim f(0)e^{C\|\theta_0\|_{L^p} t}e^{C\|\theta\|_{L^1_tB_{\infty,2}^{0,\alpha}}}.
\end{equation}
Let $N\in\NN$, then  Bernstein inequalities and Theorem \ref{max-princ} give
\begin{eqnarray*}
\nonumber\|\theta\|_{L^1_tB_{\infty,2}^{0,\alpha}}&\le&t\Big\|\Big(\sum_{q< N}(1+|q|)^{2\alpha}\|\Delta_q\theta\|_{L^\infty}^2\Big)^{\frac12}\Big\|_{L^\infty_t}+\|(\hbox{Id}-S_N)\theta\|_{L^1_t B_{\infty,1}^{0,\alpha}}\\
\nonumber&\lesssim&t\|\theta\|_{L^\infty_{t,x}}N^{\frac12+\alpha}+\sum_{q\geq N }(1+|q|)^\alpha\|\Delta_q\theta\|_{L^1_tL^\infty}\\
\nonumber&\lesssim& 
t\|\theta_0\|_{L^\infty}N^{\frac12+\alpha}+\sum_{q\geq N }(1+|q|)^\alpha\|\Delta_q\theta\|_{L^1_tL^\infty}\\
&\lesssim& N^{\frac12+\alpha}\|\theta_0\|_{L^\infty} t+\sum_{q\geq N }2^{q\frac2p}(1+|q|)^\alpha\|\Delta_q\theta\|_{L^1_tL^p}.
\end{eqnarray*}
Using Theorem \ref{thm99} then for  $p>2$    and  $0<\varepsilon<1-\frac2p$ we obtain
\begin{eqnarray*}
\nonumber
\sum_{q\geq N }(1+|q|)^\alpha2^{q\frac2p}\|\Delta_q\theta\|_{L^1_tL^p}&\lesssim&\sum_{q\geq N }(1+|q|)^{2\alpha}2^{q(\frac2p-1)}\Big(  \|\theta_0\|_{L^p}+\|\theta_0\|_{L^\infty}\|\omega\|_{L^1_tL^p}\Big)
\\
&\lesssim&\sum_{q\geq N }2^{q(\frac2p+\varepsilon-1)}\Big(  \|\theta_0\|_{L^p}+\|\theta_0\|_{L^\infty}\|\omega\|_{L^1_tL^p}\Big)\\
&\lesssim&  \|\theta_0\|_{L^p}+2^{N(-1+\varepsilon+\frac2p)}\|\theta_0\|_{L^\infty}\|\omega\|_{L^1_tL^p}.
\end{eqnarray*}
Consequently,
\begin{eqnarray}\label{es43}
\nonumber\|\theta\|_{L^1_tB_{\infty,2}^{0,\alpha}}\lesssim N^{\frac12+\alpha}\|\theta_0\|_{L^\infty} t+ \|\theta_0\|_{L^p}+2^{N(-1+\varepsilon+\frac2p)}\|\theta_0\|_{L^\infty}\|\omega\|_{L^1_tL^p}.
\end{eqnarray}
We choose $N$ as follows
 $$
 {N= \Bigg[\frac{\log\big(e+\|\omega\|_{L^1_tL^p}\big)}{(1-\varepsilon-2/p)\log2}\Bigg].}
 $$ 
 This yields
\begin{equation*}
\|\theta\|_{L^1_tB_{\infty,2}^{0,\alpha}}\lesssim\|\theta_0\|_{L^\infty\cap L^p}+\|\theta_0\|_{L^\infty} t\log^{\frac12+\alpha}\Big(e+\int_0^t\|\omega(\tau)\|_{L^p}d\tau\Big).
\end{equation*}
Combining this estimate with (\ref{gr1}) we get
\begin{eqnarray}\label{ineq23}
\nonumber\|\theta\|_{L^1_tB_{\infty,2}^{0,\alpha}}&\lesssim& \|\theta_0\|_{L^\infty\cap L^p}+\|\theta_0\|_{L^\infty} t\, \log^{\frac12+\alpha}\Big(e+Cf(0)e^{C\|\theta_0\|_{L^p}t}e^{C\|\theta\|_{L^1_tB_{\infty,2}^{0,\alpha}}}\Big)\\
&\le&{C}_0\log^{\frac12+\alpha}(e+f(0))\,(1+t^{\frac32+\alpha})+C\|\theta_0\|_{L^\infty} t\|\theta\|_{L^1_tB_{\infty,2}^{0,\alpha}}^{\frac12+\alpha},
\end{eqnarray}
where ${C}_0$ is a constant depending on $\|\theta_0\|_{L^p\cap L^\infty}$.

\underline{Case1: $\alpha<\frac12.$}
We will use the following lemma.
\begin{lem}
Let $a,b> 0$ and $\alpha\in[0,1[$. Let $x\in\RR_+$ be a solution of the inequality
$$
(\star)\quad x\le a+b x^\alpha.
$$
Then there exists $C:=C(\alpha)$ such that 
$$
x\le C_\alpha(a+b^{\frac{1}{1-\alpha}}).
$$
\end{lem}
\begin{proof}
We set $y=a^{-1}x$. Then the inequality $(\star)$ becomes
$$
y\le 1+ba^{\alpha-1}y^\alpha.
$$ 
We will look for a number $\mu>0$ such that $y\le e^\mu.$ Then it suffices to find $\mu$ such that
$$
1+ba^{\alpha-1}e^{\mu\alpha}\le e^\mu.
$$
It suffices also to find $\mu$ such that
$$
(1+ba^{\alpha-1})e^{\mu\alpha}= e^\mu.
$$
This gives $e^{\mu}=(1+ba^{\alpha-1})^{\frac{1}{1-\alpha}}$. Now we can use the inequality: for every $t,s\geq0$
$$
(t+s)^{\frac{1}{1-\alpha}}\le C_\alpha(t^{\frac{1}{1-\alpha}}+s^{\frac{1}{1-\alpha}}).
$$
It follows that
$$
y\le C_\alpha(1+a^{-1}b^{\frac{1}{1-\alpha}}).
$$
This yields 
$$
x\le C_\alpha (a+b^{\frac{1}{1-\alpha}})
$$
\end{proof}
Applying this lemma to (\ref{ineq23}) we get for every $t\in\RR_+$
\begin{eqnarray}\label{eq34}
\nonumber\|\theta\|_{L^1_tB_{\infty,2}^{0,\alpha}}&\le& C_0 (t^\frac32+t^{\frac{2}{1-2\alpha}})\\
\nonumber&\le&C_0 (1+t^{\frac{2}{1-2\alpha}})\\
&\le&\Phi_1(t).
\end{eqnarray}
It follows from (\ref{gr1})
\begin{eqnarray}
\label{gr6}
\nonumber f(t)&\le& C_0e^{C_0t^{\frac{2}{1-2\alpha}}}\\
&\le&\Phi_2(t)
\end{eqnarray}
Applying Theorem \ref{thm99} and (\ref{gr6}) we  get  for every $\epsilon>0, q\in\NN$
\begin{eqnarray}
\label{gr4}
\nonumber2^{q}(1+|q|)^{-\alpha}\|\Delta_q\theta\|_{L^1_tL^p}&\leq& C_0e^{C_0t^{\frac{2}{1-2\alpha}}}\\
&\le&\Phi_2(t).
\end{eqnarray}
\underline{Case 2: $\alpha=\frac12$.} The estimate (\ref{ineq23}) becomes
$$
\|\theta\|_{L^1_tB_{\infty,2}^{0,\frac12}}\le {C}_0\log(e+f(0))(1+t^{2})+C\|\theta_0\|_{L^\infty} t\|\theta\|_{L^1_tB_{\infty,2}^{0,\frac12}},
$$
with ${C}_0$  a constant depending on $\|\theta_0\|_{L^p\cap L^\infty}$.
Hence if we choose $t$ small enough such that 
\begin{equation}\label{eq35}
C\|\theta_0\|_{L^\infty} t=\frac12,
\end{equation}
then
$$
\|\theta\|_{L^1_tB_{\infty,2}^{0,\alpha}}\le {C}_0\log(e+f(0)).
$$
From \eqref{gr1} we get that 
$$
f(t)\le{C}_0(e+f(0))^{C_0}.
$$
Now let $t$ be a given positive time and choose a partition $(t_i)_{i=1}^N$ of $[0,t]$ such that
\begin{equation}\label{parti}
C\|\theta_0\|_{L^\infty} (t_{i+1}-t_i)\approx\frac12.
\end{equation}
Set $a_i:=\int_{t_i}^{t_{i+1}}\|\theta(\tau)\|_{B_{\infty,2}^{0,\frac12}}d\tau$ and $b_i=f(t_i)$.  Thus reproducing  similar computations to   (\ref{ineq23}) yields
\begin{eqnarray*}
a_i\le{C}_0\log(e+b_i)(1+(t_{i+1}-t_i)^{2})+C\|\theta_0\|_{L^\infty}(t_{i+1}-t_i)a_i.
\end{eqnarray*}
 Hence we get
\begin{equation}\label{recc1}
a_i\le {C}_0\log(e+b_i).
\end{equation}
The analogous  estimate to  (\ref{gr1}) is
\begin{eqnarray}\label{recc2}
\nonumber b_{i+1}&\lesssim& b_i e^{C(t_{i+1}-t_i)\|\theta_0\|_{L^p}}e^{Ca_i}\\
&\le&C_0 b_i e^{Ca_i}.
\end{eqnarray}
Combining (\ref{recc1}) and (\ref{recc2}) yields
$$
b_{i+1}\le C_0(e+b_i)^{C_0}.
$$
By induction we can prove that for every $i\in\{1,..,N\}$ we have
$$
b_i\le C_0e^{\exp{C_0 i}}
$$
and consequently from \eqref{recc1}
$$
a_i\le C_0e^{C_0 i}.
$$
It follows that
\begin{eqnarray*}
\|\theta\|_{L^1_tB_{\infty,2}^{0,\frac12}}&=&\sum_{i=^1}^Na_i\\
&\le&C_0e^{C_0 N}\\
&\leq& C_0e^{C_0 t}.
\end{eqnarray*}
We have used in the last inequality the fact that
$$
N\approx C_0 t
$$
which is a consequence of (\ref{parti}). We have also obtained
$$
f(t)\le C_0e^{\exp{C_0 t}}.
$$
It is not hard to see  that from \eqref{gr4} one can obtain that  for every $s<1$
\begin{equation}
\label{x5}
\|\theta\|_{ L^1_tB_{p,1}^s}\leq  \|\theta\|_{ \widetilde{ L}^1_tB_{p,\infty}^{1,-\alpha}} \leq  \Phi_2(t).
\end{equation}

This ends the proof of Proposition \ref{max-pro}.
\begin{Rema} Combining \eqref{x5}  with  Bernstein inequalities and the fact that $p>2$ this yields
\begin{equation}
\label{x9}
\|\theta\|_{ L^1_tB_{\infty,1}^\epsilon}\leq
 \Phi_2(t),
\end{equation}
for every $\epsilon <1-\frac2p$.
\end{Rema}
\vspace{0.5cm}
\end{proof}

\subsubsection{$L^\infty$-bound of the vorticity.} We will prove the following result.
 
\begin{prop}
\label{pr0} 
Let $\alpha\in[0,\frac12], \lambda\geq e^{3+2\alpha}, p\in]2,\infty[$ and  $(v,\theta)$ be a smooth solution of  the system  \eqref{Bouss}  such that   $\omega^0,\theta_0,\mathcal{R}_\alpha\theta_0\in L^p\cap L^\infty$.  Then
we have
\begin{equation}
\label{x10}
\|\omega(t)\|_{L^\infty}+\|\mathcal{R}_\alpha\theta(t)\|_{L^\infty}\le \Phi_3(t)
\end{equation}
and 
\begin{equation}
\label{x11}
\|v(t)\|_{L^\infty}\le \Phi_4(t).
\end{equation}
\end{prop}
\begin{proof}

\
 {\sl\underline{ Proof of \eqref{x10}.}}
By using the  maximum principle for the transport equation \eqref{bm1},  we get
$$
\|\Gamma(t)\|_{L^\infty}\leq\|\Gamma^0\|_{L^\infty}+\int_0^t\|[\mathcal{R}_\alpha, v\cdot\nabla]\theta(\tau)\|_{L^\infty}
d\tau.
$$
Since the function $\mathcal{R}_\alpha\theta$ satisfies the equation
\begin{equation*}
\label{rtheta}
\big(\partial_t+v\cdot\nabla+\DD\log^{-\alpha}(\lambda+\DD)\big)\mathcal{R}_\alpha\theta=-[\mathcal{R}_\alpha,v\cdot\nabla]\theta,
\end{equation*}
then using Theorem \ref{max-princ}  we get
$$
\|\mathcal{R}_\alpha\theta(t)\|_{L^\infty}\leq\|\mathcal{R}_\alpha\theta(t)\|_{L^\infty}+\int_0^t\|[\mathcal{R}_\alpha, v\cdot\nabla]\theta(\tau)\|_{L^\infty}
d\tau.
$$
Thus we obtain
\begin{eqnarray*}
\|\Gamma(t)\|_{L^\infty}+\|\mathcal{R}_\alpha\theta(t)\|_{L^\infty}&\leq&\|\Gamma^0\|_{L^\infty}+\|\mathcal{R}_\alpha\theta_0\|_{L^\infty}+2\int_0^t\|[\mathcal{R}_\alpha, v\cdot\nabla]\theta(\tau)\|_{L^\infty}
d\tau
\\
&\leq & C_0+\int_0^t  \|[\mathcal{R}_\alpha, v\cdot\nabla]\theta(\tau)\|_{B_{\infty,1}^0}
d\tau.
\end{eqnarray*}
It follows from  Theorem \ref{max-princ}, Proposition \ref{propcom}-(2)  and  Proposition \ref{max-pro}   
\begin{eqnarray*}
\|\omega(t)\|_{L^\infty}+\|\mathcal{R}_\alpha\theta(t)\|_{L^\infty}&\lesssim&C_0+\int_0^t\|\omega(\tau)\|_{L^\infty\cap L^p}\big(\|\theta(\tau)\|_{B_{\infty,1}^\epsilon}+\|\theta(\tau)\|_{L^p}\big)
d\tau
\\
&\lesssim&C_0+\|\omega\|_{L^\infty_t L^p}\big(\|\theta\|_{L^1_tB_{\infty,1}^\epsilon}+t\|\theta_0\|_{L^p}\big)\\
&+&\int_0^t\|\omega(\tau)\|_{L^\infty}\big(\|\theta(\tau)\|_{B_{\infty,1}^\epsilon}+\|\theta_0\|_{L^p}\big)
d\tau.
\end{eqnarray*}
Let $0<\epsilon<1-\frac2p$ then using \eqref{x9} we get
$$
\|\omega(t)\|_{L^\infty}+\|\mathcal{R}_\alpha\theta(t)\|_{L^\infty}\lesssim\Phi_2(t)+\int_0^t\|\omega(\tau)\|_{L^\infty}\big(\|\theta(\tau)\|_{B_{\infty,1}^\epsilon}+\|\theta_0\|_{L^p}\big)d\tau.
$$
 Therefore we obtain by    Gronwall lemma  and  \eqref{x9}  
 \begin{eqnarray*}
\|\omega(t)\|_{L^\infty}+\|\mathcal{R}_\alpha\theta(t)\|_{L^\infty}\leq \Phi_3(t).
\end{eqnarray*}

{\sl \underline{Proof of \eqref{x11}.}}
Let  $N\in\NN$ to be chosen later.  Using the fact that $\|\dot\Delta_qv\|_{L^\infty}\approx 2^{-q}\|\dot\Delta_q\omega\|_{L^\infty}$, then we get\begin{eqnarray*}
\|v(t)\|_{L^\infty}&\leq& \|\chi(2^{N}\DD)v(t)\|_{L^\infty}+\sum_{q\geq-N}2^{-q}\|\dot\Delta_q\omega(t)\|_{L^\infty}\\
\\
&\le&  \|\chi(2^{N}\DD)v(t)\|_{L^\infty}+2^N\|\omega(t)\|_{L^\infty}.
\end{eqnarray*} 
Applying the frequency localizing operator to the velocity equation we get
$$
\chi(2^{N}\DD)v=\chi(2^{N}\DD)v_0+\int_0^t\mathcal P\chi(2^{N}\DD)\theta(\tau)d\tau+\int_0^t\mathcal P\chi(2^{N}\DD){\rm div}(v\otimes v)(\tau)d\tau,
$$
where $\mathcal P$ stands for Leray projector. From   Lemma \ref{lb}, Calder\'on-Zygmund estimate and the uniform boundness of $\chi(2^{N}\DD)$  we get 
\begin{eqnarray*}
\int_0^t\|\chi(2^{N}\DD)\mathcal P\theta(\tau)\|_{L^\infty}d\tau&\lesssim&2^{-N\frac2p} \int_0^t\|\theta(\tau)\|_{L^p}d\tau\\
&\lesssim&t\|\theta_0\|_{L^p}.
\end{eqnarray*}
Using Proposition \ref{cor1}-(2) we find
$$
\int_0^t\|\mathcal P\chi(2^{N}\DD){\rm div}(v\otimes v)(\tau)\|_{L^\infty}d\tau\lesssim 2^{N}\int_0^t\|v(\tau)\|_{L^\infty}^2d\tau.
$$
The outcome is 
\begin{eqnarray*}
\|v(t)\|_{L^\infty}&\lesssim& \|v_0\|_{L^\infty}+t \| \theta_0\|_{L^p}+ 2^{-N}\int_0^t \| v(\tau)\|_{L^\infty}^2d\tau+2^N\|\omega(t)\|_{L^\infty}
\\
&\lesssim&   2^{-N}\int_0^t \| v(\tau)\|_{L^\infty}^2d\tau+  2^N  \Phi_3(t) 
\end{eqnarray*} 
Choosing judiciously $N$ we find
$$
\|v(t)\|_{L^\infty}\le\Phi_3(t)\Big(1+\Big(\int_0^t\|v(\tau)\|_{L^\infty}^2d\tau\Big)^{\frac12}\Big).
$$From Gronwall lemma we get
$$
\|v(t)\|_{L^\infty}\le \Phi_4(t).
$$
\end{proof}
\subsubsection{Lipschitz bound of the velocity}
Now we will establish the following result.
\begin{prop}
\label{pr10}
Let $\alpha\in[0,\frac12],\lambda\geq e^{3+2\alpha}, p\in]2,\infty[$ and 
$(v,\theta)$ be a smooth solution of the system \eqref{Bouss} with $\omega^0,\theta_0,\mathcal{R}_\alpha\theta_0\in B_{\infty,1}^0\cap L^p$. Then
$$
\|\mathcal{R}_\alpha\theta(t)\|_{B_{\infty,1}^0}+\|\omega(t)\|_{B_{\infty,1}^0}+\|v(t)\|_{B_{\infty,1}^1}\leq \Phi_4(t).
$$

\end{prop}
\begin{proof} 
Applying 
Corollary \ref{thmlog} to the equations \eqref{bm} and \eqref{bm},  we obtain
\begin{equation}\label{eqlog}
\|\Gamma(t)\|_{B_{\infty,1}^0}+\|\mathcal R_\alpha\theta(t)\|_{B_{\infty,1}^0} \lesssim\Big(C_0+\big\|[\mathcal{R}_\alpha, v\cdot\nabla]\theta\big\|_{L^1_tB_{\infty,1}^0}\Big)\Big(1+\|\nabla v\|_{L^1_tL^\infty}  \Big).
\end{equation}
Thanks to  Theorem \ref{propcom}, Propositions \ref{pr0}, \ref{max-pro} and \eqref{x9} we get
\begin{eqnarray*}
\big\|[\mathcal{R}_\alpha, v\cdot\nabla]\theta\big\|_{L^1_tB_{\infty,1}^0}&\lesssim &\int_0^t(\|\omega(\tau)\|_{L^\infty}+\|\omega(\tau)\|_{L^p})\big(\|\theta(\tau)\|_{B_{\infty,1}^\epsilon}+\|\theta(\tau)\|_{L^p}\big)d\tau
\\
&\lesssim & \Phi_3(t).
\end{eqnarray*}
By easy computations  we get
\begin{eqnarray}
\nonumber
\|\nabla v\|_{L^\infty}&\le& \|\nabla\Delta_{-1}v\|_{L^\infty}+\sum_{q\in\NN}\|\Delta_q\nabla v\|_{L^\infty}
\\
\nonumber
&\lesssim& \|\omega\|_{L^p}+\sum_{q\in\NN}\|\Delta_q\omega\|_{L^\infty}
\\
\label{x8}
&\lesssim&  \Phi_2(t)+\|\omega(t)\|_{B_{\infty,1}^0}.
\end{eqnarray}
Putting together \eqref{eqlog} and \eqref{x8} leads to
\begin{eqnarray*}
\|\omega(t)\|_{B_{\infty,1}^0}\leq \|\Gamma(t)\|_{B_{\infty,1}^0}+\|\mathcal R_\alpha\theta(t)\|_{B_{\infty,1}^0}\leq  \Phi_3(t)\Big(1+\int_0^t\|\omega(\tau)\|_{B_{\infty,1}^0}d\tau\Big).
\end{eqnarray*}
Thus we obtain from Gronwall inequality
\begin{equation}\label{x12}
\|\omega(t)\|_{B_{\infty,1}^0}+\|\mathcal R_\alpha\theta(t)\|_{B_{\infty,1}^0}\leq\Phi_4(t).
\end{equation}
Coming back to \eqref{x8} we get
\begin{eqnarray*}
\nonumber
\|\nabla v(t)\|_{L^\infty}\leq   \Phi_4(t).
\end{eqnarray*}
Let us move to the estimate of $v$ in the space $B_{\infty,1}^1$. By definition we have
$$
\|v(t)\|_{B_{\infty,1}^1}\lesssim \|v(t)\|_{L^\infty}+\|\omega(t)\|_{B_{\infty,1}^0}.
$$
Combined with \eqref{x11} and   \eqref{x12} this yields 
$$
\|v(t)\|_{B_{\infty,1}^1}\leq \Phi_4(t)
$$
The proof of \mbox{Proposition \ref{pr10}} is now achieved.  
\end{proof}

\subsection{Uniqueness}
We will   show that the Boussinesq system \eqref{Bouss} has a unique solution in the following function space
$$\mathcal{E}_T= (L^\infty_TB_{\infty,1}^0\cap L^1_TB_{\infty,1}^1)\times (L^\infty_TL^p\cap \widetilde L^1_T B_{p,\infty}^{1,-\alpha}),\quad 2<p<\infty.
$$
Let $(v^1,\theta^1)$ and  $(v^2,\theta^2)$  be two solutions of \eqref{Bouss} belonging to the space $\mathcal{E}_T$ and denote 
$$v=v^2-v^1,\quad\theta=\theta^2-\theta^1.
$$ Then we get 
\begin{equation*}
\left\{ 
\begin{array}{ll}
 \partial_t v+v^2\cdot\nabla v=-\nabla \pi-v\cdot\nabla v^1+\theta e_2\\ 
\partial_t\theta+v^2\cdot\nabla \theta+\frac{|\textnormal{D}|}{\log^\alpha(\lambda+\DD)}\theta=-v\cdot\nabla\theta^1\\
v_{| t=0}=v_0, \quad \theta_{| t=0}=\theta_0. 
\end{array} \right.
\end{equation*}
According to  Proposition  \ref{prop-Bes} we have
$$
\|v(t)\|_{B_{\infty,1}^0}\le Ce^{CV_1(t)}\Big(\|v_0\|_{B_{\infty,1}^0}+\|\nabla \pi\|_{L^1_tB_{\infty,1}^0}+\| v\cdot\nabla v^1\|_{L^1_tB_{\infty,1}^0}+\|\theta\|_{L^1_tB_{\infty,1}^0}\Big),
$$
with $V_1(t)=\|\nabla v^1\|_{L^1_tL^\infty}.$
Straightforward computations using the incompressibility of the flows gives
\begin{eqnarray*}
\nabla \pi&=&-\nabla\Delta^{-1}\textnormal{div }(v\cdot\nabla(v^1+v^2))+\nabla\Delta^{-1}\partial_2\theta\\
&=&\hbox{I}+\hbox{II}.
\end{eqnarray*}
To estimate the first term of the RHS we use the definition
$$
\|\hbox{I}\|_{B_{\infty,1}^0}\lesssim\|(\nabla\Delta^{-1}\textnormal{div })\textnormal{div }\Delta_{-1}(v\otimes (v^1+v^2))\|_{L^\infty}+\|v\cdot\nabla(v^1+v^2)\|_{B_{\infty,1}^1}
$$
From Proposition 3.1-(2)  of \cite{HKR} and Besov embeddings  we have
\begin{eqnarray*}
\|(\nabla\Delta^{-1}\textnormal{div })\textnormal{div }\Delta_{-1}(v\otimes (v^1+v^2))\|_{L^\infty}&\lesssim& \|v\otimes (v^1+v^2)\|_{L^\infty}\\
&\lesssim&\|v\|_{B_{\infty,1}^0} \|v^1+v^2\|_{B_{\infty,1}^0}.
\end{eqnarray*}
Using the incompressibility of $v$ and using Bony's decomposition one can easily obtain
$$
\|v\cdot\nabla(v^1+v^2)\|_{B_{\infty,1}^0}\lesssim \|v\|_{B_{\infty,1}^0}\|v^1+v^2\|_{B_{\infty,1}^1}.
$$
Putting together  these estimates yields
\begin{equation}\label{s12}
\|\hbox{I}\|_{B_{\infty,1}^0}\lesssim \|v\|_{B_{\infty,1}^0}\|v^1+v^2\|_{B_{\infty,1}^1}.
\end{equation}
Let us now show how to estimate the second term $\hbox{II}$. By using Besov embeddings and 
 Calder\'on-Zygmund estimate we get
 \begin{eqnarray*}
 \|\hbox{II}\|_{B_{\infty,1}^0}&\lesssim& \|\nabla\Delta^{-1}\partial_2\theta\|_{B_{p,1}^{\frac2p}}\\
 &\lesssim&  \|\theta\|_{B_{p,1}^{\frac2p}}.
 \end{eqnarray*}
 Combining this estimate with (\ref{s12}) yields
 \begin{eqnarray}
 \label{1170}
\nonumber\|v(t)\|_{B_{\infty,1}^0}\lesssim e^{CV(t)}\Big(\|v_0\|_{B_{\infty,1}^0}&+&\int_0^t\| v(\tau)\|_{B_{\infty,1}^0}\big[1+\|(v^1,v^2)(\tau)\|_{B_{\infty,1}^1}\big]d\tau\Big)\\&+&e^{CV(t)}
\|\theta\|_{L^1_tB_{p,1}^{\frac2p}},
 \end{eqnarray}
where $V(t):=\|(v^1,v^2)\|_{L^1_tB_{\infty,1}^1}.$

Now we have  to estimate $\|\theta\|_{L^1_tB_{p,1}^{\frac2p}}$.  By applying $\Delta_q$ to the equation of $\theta $ and arguing similarly to the proof of Theorem \ref{thm99} we obtain for $q\in\NN$
\begin{eqnarray*}
\nonumber\|\theta_q(t)\|_{L^p} &\lesssim&  e^{-ct2^q(1+q)^{-\alpha}}\|\theta_q^0\|_{L^p}+\int_0^t e^{-c2^q(1+q)^{-\alpha}(t-\tau)} \|\Delta_q(v\cdot\nabla \theta^1)(\tau)\|_{L^p}d\tau\\
&+&\int_0^te^{-c2^q(1+q)^{-\alpha}(t-\tau)} \|\big[v^2\cdot\nabla,\Delta_q\big]\theta(\tau) \|_{L^p}d\tau.
\end{eqnarray*}
Remark, first, that an obvious H\"older inequality yields that for every $\varepsilon\in[0,1]$ there exists an absolute constant $C$ such that 
$$
\int_0^t e^{-c\tau2^q(1+q)^{-\alpha}}d\tau\leq C t^{1-\varepsilon}2^{-q\varepsilon}(1+q)^{\alpha\,\varepsilon},\qquad \forall\, t\geq 0.
$$
Using this fact  and integrating in time 
\begin{eqnarray}
\label{113}
\nonumber2^{q\frac2p}\|\theta_q\|_{L^1_tL^p}& \lesssim &  (q+1)^{\alpha}2^{q(-1+\frac2p)}\|\theta_q^0\|_{L^p}\\
\nonumber&+& t^{1-\varepsilon} (q+1)^{\alpha\varepsilon}2^{q(-\varepsilon+\frac2p)}\int_0^t \Big( \|\Delta_q(v\cdot\nabla \theta^1)(\tau)\|_{L^p}+ \|[v^2\cdot\nabla,\Delta_q]\theta(\tau) \|_{L^p}\Big)d\tau\\
&=& (q+1)^{\alpha}2^{q(-1+\frac2p)}\|\theta_q^0\|_{L^p}+{\rm I}_q(t)+{\rm II}_q(t).
\end{eqnarray}
Using Bony's decomposition we get easily
\begin{eqnarray*}
 \|\Delta_q(v\cdot\nabla \theta^1)(t)\|_{L^p}& \lesssim&   \|v(t)\|_{L^\infty}\sum_{j\leq q+2} 2^{j}  \|\Delta_j\theta^1(t)\|_{L^p}\\
 &+& 2^q \|v(t)\|_{L^\infty}\sum_{j\geq q-4}   \|\Delta_j\theta^1(t)\|_{L^p}\\
 &\lesssim&\|v(t)\|_{L^\infty}\sum_{j\leq q+2}(1+|j|)^\alpha \big(2^{j}(1+|j|)^{-\alpha}  \|\Delta_j\theta^1(t)\|_{L^p}\big)\\
 &+&  \|v(t)\|_{L^\infty}\sum_{j\geq q-4} 2^{q-j}(1+|j|)^{\alpha} \big(2^j(1+|j|)^{-\alpha} \|\Delta_j\theta^1(t)\|_{L^p}\big).
  \end{eqnarray*}
  Integrating in time we get
  \begin{eqnarray}\label{114}
 \nonumber {\rm{I}}_q (t) & \lesssim& t^{1-\varepsilon} \|v\|_{L^\infty_tL^\infty}2^{q(\frac2p-\varepsilon)}(q+1)^{1+\alpha(1+\varepsilon)}\|\theta^1\|_{\widetilde L^1_tB_{p,\infty}^{1,-\alpha}}\\ 
 \nonumber&+& t^{1-\varepsilon} \|v\|_{L^\infty_tL^\infty}\|\theta^1\|_{\widetilde L^1_tB_{p,\infty}^{1,-\alpha}}2^{q(\frac2p+1-\varepsilon)}(q+1)^{\alpha(1+\varepsilon)}\sum_{j\geq q-4}2^{-j} (1+|j|)^{\alpha} \\
 &\lesssim&t^{1-\varepsilon} \|v\|_{L^\infty_tL^\infty}2^{q(\frac2p-\varepsilon)}(q+1)^{1+\alpha(1+\varepsilon)}\|\theta^1\|_{\widetilde L^1_tB_{p,\infty}^{1,-\alpha}}.
  \end{eqnarray}
  To estimate the term ${\rm II}_q$ we use the following classical commutator ( since $2/p<1$), see \cite{che1}
  $$
  \|[v^2\cdot\nabla,\Delta_q]\theta \|_{L^p}\lesssim  2^{-q\frac2p}\|\nabla v^2\|_{L^\infty}\|\theta\|_{B_{p,1}^{\frac2p}}.
  $$
  
  Thus we obtain,
\begin{equation}\label{115}
{\rm II}_q (t) \lesssim t^{1-\varepsilon}(q+1)^{\alpha\varepsilon} 2^{-q\varepsilon}\|\nabla v^2\|_{L^\infty_tL^\infty}\|\theta\|_{L^1_tB_{p,1}^{\frac2p}}.
\end{equation}

We choose $\varepsilon\in]0,1[$ such that $\frac2p-\varepsilon<0,$ which is possible since $p>2.$ Then  combining  (\ref{113}),  (\ref{114}) and (\ref{115}) we get 
\begin{eqnarray*}
\|\theta\|_{L^1_tB_{p,1}^{\frac2p}}&\lesssim & \|\theta_0\|_{L^p}+t^{1-\varepsilon} \|v\|_{L^\infty_tL^\infty}\|\theta^1\|_{\widetilde L^1_tB_{p,\infty}^{1,-\alpha}}+t^{1-\varepsilon}\|\nabla v^2\|_{L^\infty_tL^\infty}\|\theta\|_{L^1_tB_{p,1}^{\frac2p}}.
\end{eqnarray*}
It follows that there exists small $\delta>0$ such that  for $t\in[0,\delta]$ 
$$
\|\theta\|_{L^1_tB_{p,1}^{\frac2p}}\lesssim  \|\theta_0\|_{L^p}+t^{1-\varepsilon} \|v\|_{L^\infty_tL^\infty}\|\theta^1\|_{\widetilde L^1_tB_{p,\infty}^{1,-\alpha}}.
$$ 
Plugging this estimate into \eqref{1170} we find
$$
\|v\|_{L^\infty_tB_{\infty,1}^0}\lesssim e^{CV(t)}\Big(\|v_0\|_{B_{\infty,1}^0}+\|\theta_0\|_{L^p}+t\| v\|_{L^\infty_tB_{\infty,1}^0}+ t^{\varepsilon} \|v\|_{L^\infty_tL^\infty}\|\theta^1\|_{\widetilde L^1_tB_{p,\infty}^{1,-\alpha}}\Big).
$$
If $\delta$ is sufficiently small then we get for $t\in[0,\delta]$
\begin{equation}\label{uni1}
\|v\|_{L^\infty_tB_{\infty,1}^0}\lesssim \|v_0\|_{B_{\infty,1}^0}+\|\theta_0\|_{L^p}.
\end{equation}
This gives in turn
\begin{equation}\label{uni2}
\|\theta\|_{L^1_tB_{p,1}^{\frac2p}}\lesssim \|v_0\|_{B_{\infty,1}^0}+\|\theta_0\|_{L^p}.
\end{equation}
This gives in particular the uniqueness on  $[0,\delta]$.  Iterating this argument  yields   the uniqueness in $[0,T]$. 

\subsection{Existence}
We consider the following system
\begin{equation} 
\left\{ \begin{array}{ll} 
\partial_{t}v_n+v_n\cdot\nabla v_n+\nabla \pi_n=\theta_n e_{2}\\ 
\partial_{t}\theta_n+v_n\cdot\nabla\theta_n+\frac{\vert \textnormal{D}}{\log^\alpha(\lambda+\DD)}\theta_n=0\\
\textnormal{div}v_n=0\\
{v_n}_{| t=0}=S_nv^{0}, \quad {\theta_n}_{| t=0}=S_n\theta^{0}. 
\end{array} \right. \tag{B$_{n}$}
\end{equation}
By using the same method as \cite{hk} we can prove that this system has a unique  local smooth solution $(v_n,\theta_n)$. The global existence of these solutions is governed by  the following criterion: we can push the construction beyond the time $T$ if  the quantity $\|\nabla v_n\|_{L^1_T L^\infty}$ is finite. Now from the {\it a priori} estimates the Lipschitz norm can not blow up in finite time and then the solution $(v_n,\theta_n)$ is globally defined. Once again from the {\it a priori} estimates   we have  for $2<p<\infty$ 
$$
\|v_n\|_{ L^\infty_TB_{\infty,1}^1}+\|\omega_n\|_{L^\infty_T L^p}+\|\theta_n\|_{L^\infty_T\mathcal{X}_p}\leq \Phi_4(T).
$$
The space $\mathcal{X}_p$ was introduced before the statement of Theorem \ref{theo1}
It follows that up to an extraction the sequence $(v_n,\theta_n)$ is weakly convergent to  $(v,\theta)$ belonging to $L^\infty_TB_{\infty,1}^1\times  L^\infty_T\mathcal{X}_p,$ with $\omega\in L^\infty_T L^p$.  For $(n,m)\in\NN^2$ we set $v_{n,m}=v_n-v_m$ and $\theta_{n,m}=\theta_n-\theta_m$ then  according to the estimate (\ref{uni1}) and (\ref{uni2}) we get  for $T=\delta$
$$\|v_{n,m}\|_{L^\infty_TB_{\infty,1}^0}+\|\theta_{n,m}\|_{L^1_TB_{p,1}^{\frac2p}}\lesssim \|S_nv_0-S_m v_0\|_{B_{\infty,1}^0}+\|S_n\theta_0-S_m \theta_0\|_{L^p}.
$$
This shows that $(v_n,\theta_n)$ is a Cauchy sequence in the Banach space $L^\infty_TB_{\infty,1}^0\times L^1_TB_{p,1}^{\frac2p}$ and then it converges strongly to $(v,\theta).$ This allows to pass to the limit in the system $({\rm B}_n)$ and then we get that $(v,\theta)$ is a solution of the Boussinesq system (\ref{Bouss}).

\end{document}